\documentclass{amsart}%
\usepackage{amssymb}
\usepackage{amsmath}
\usepackage{amsfonts}
\usepackage{graphicx}%
\providecommand{\U}[1]{\protect\rule{.1in}{.1in}}
\providecommand{\U}[1]{\protect \rule{.1in}{.1in}}
\providecommand{\U}[1]{\protect \rule{.1in}{.1in}}
\providecommand{\U}[1]{\protect \rule{.1in}{.1in}}
\providecommand{\U}[1]{\protect \rule{.1in}{.1in}}
\providecommand{\U}[1]{\protect \rule{.1in}{.1in}}
\providecommand{\U}[1]{\protect \rule{.1in}{.1in}}
\providecommand{\U}[1]{\protect \rule{.1in}{.1in}}
\providecommand{\U}[1]{\protect \rule{.1in}{.1in}}
\providecommand{\U}[1]{\protect \rule{.1in}{.1in}}
\providecommand{\U}[1]{\protect \rule{.1in}{.1in}}
\providecommand{\U}[1]{\protect \rule{.1in}{.1in}}
\providecommand{\U}[1]{\protect \rule{.1in}{.1in}}
\providecommand{\U}[1]{\protect \rule{.1in}{.1in}}
\providecommand{\U}[1]{\protect \rule{.1in}{.1in}}
\providecommand{\U}[1]{\protect \rule{.1in}{.1in}}

\theoremstyle{plain}

\numberwithin{equation}{section}
\begin{document}
\title[Inverse Problems For Dirac Operators ]{Inverse Problems For Dirac Operators With a Finite Number of Transmission
Conditions }
\author{Yal\c{c}\i n G\"{U}LD\"{U} and Merve ARSLANTA\c{S} }
\address[Yal\c{c}\i n G\"{U}LD\"{U}]{, [Merve ARSLANTA\c{S}] Department of Mathematics,
Faculty of Science, Cumhuriyet University, 58140 Sivas, Turkey}
\email{yguldu@gmail.com, merveguray0@gmail.com}
\subjclass[2000]{ Primary 34A55,{\ Secondary 34B24, 34L05}}
\keywords{Dirac operator, Eigenvalues, Eigenfunctions, Transmission Conditions, Weyl Function}

\begin{abstract}
In this paper, we consider a discontinuous Dirac operator depending
polynomially on the spectral parameter and a finite number of transmission
conditions. We get some properties of eigenvalues and eigenfunctions. Then, we
investigate some uniqueness theorems by using Weyl function and some spectral data.

\end{abstract}
\maketitle

\section{\textbf{Introduction}}

Inverse problems of spectral analysis recover operators by their spectral
data. Fundamental and vast studies about the classical Sturm-Liouville, Dirac
operators, Schr\"{o}dinger equation and hyperbolic equations are well studied
(see [1-7] and references therein).

Studies where eigenvalue dependent appears not only in the differential
equation but also in the boundary conditions have increased in recent years
(see [8-16] and corresponding bibliography cited therein). Moreover, boundary
conditions which depend linearly and nonlinearly on the spectral parameter are
considered in [8,16-20] and [21-27] respectively \ Furthermore, boundary value
problems with transmission conditions are also increasingly studied. These
types of studies introduce qualitative changes in the exploration. Direct and
inverse problems for Sturm-Liouville and Dirac operators with transmission
conditions are investigated in some papers (see [7, 28--31] and references
therein). Then, differential equations with the spectral parameter and
transmission conditions arise in heat, mechanics, mass transfer problems, in
diffraction problems and in various physical transfer problems (see [18, 28,
32-39] and corresponding bibliography).

More recently, some boundary value problems with eigenparameter in boundary
and transmission conditions are spread out to the case of two, more than two
or a finite number of transmission in [40-44] and references therein.

The present paper deals with the discontinuous Dirac operator depending
polynomially on the spectral parameter and a finite number of transmission
conditions. The aim of the present paper is to obtain the asymptotic formulae
of eigenvalues, eigenfunctions and to prove some uniqueness theorems.
Especially, some parameters of the considered problem can be determined by
Weyl function and some spectral data.\pagebreak

We consider a discontinuous boundary value problem $L$ with function $\rho
(x)$;%
\begin{equation}
ly:=By^{\prime}(x)+\Omega(x)y(x)=\lambda\rho(x)y(x),\ \ x\in\cup_{i=0}%
^{n}\left(  \xi_{i},\xi_{i+1}\right)  =\Lambda,\text{ }\xi_{0}=a,\text{ }%
\xi_{n+1}=b \tag{1}%
\end{equation}
where $\rho(x)=\left\{
\begin{array}
[c]{c}%
\rho_{0}\text{, }a\leq x<\xi_{1}\\
\rho_{i}\text{, }\xi_{i}<x<\xi_{i+1}\\
\rho_{n}\text{, }\xi_{n}<x\leq b
\end{array}
\right.  ,$ $i=\overline{1,n-1}\smallskip$ and $\rho_{i},$ $i=\overline{0,n}$
are given positive real numbers; $\Omega(x)=\left(
\begin{array}
[c]{ll}%
p\left(  x\right)  & q(x)\\
q(x) & r\left(  x\right)
\end{array}
\right)  ,p\left(  x\right)  ,q(x),r\left(  x\right)  \in L_{2}\left[
\Lambda,%
\mathbb{R}
\right]  ;$ $\ $\newline$B=\left(
\begin{array}
[c]{ll}%
0 & 1\\
-1 & 0
\end{array}
\right)  ,$ $y\left(  x\right)  =\left(
\begin{array}
[c]{l}%
y_{1}\left(  x\right) \\
y_{2}\left(  x\right)
\end{array}
\right)  $, $\lambda\in%
\mathbb{C}
$ is a real spectral parameter; boundary conditions$\smallskip$ at the
endpoints%
\begin{equation}
l_{1}y:=a_{2}(\lambda)y_{2}(a)-a_{1}(\lambda)y_{1}(a)=0 \tag{2}%
\end{equation}%
\begin{equation}
l_{2}y:=b_{2}(\lambda)y_{2}(b)-b_{1}(\lambda)y_{1}(b)=0 \tag{3}%
\end{equation}
\newline with transmission conditions at $n$ points $x=\xi i,$ $i=\overline
{1,n}$%
\begin{align}
l_{3}y  &  :=y_{1}(\xi i+0)-\theta_{i}y_{1}(\xi_{i}-0)=0\smallskip\tag{4}\\
l_{4}y  &  :=y_{2}(\xi_{i}+0)-\theta_{i}^{-1}y_{2}(\xi_{i}-0)-\gamma
_{i}\left(  \lambda\right)  y_{1}(\xi_{i}-0)=0\smallskip\tag{5}%
\end{align}
where $\theta_{i}$, $\xi i$ $\left(  i=\overline{1,n}\right)  $ are real
numbers, $a_{i}(\lambda),$ $b_{i}(\lambda),\left(  i=1,2\right)  $ and
$\gamma_{i}\left(  \lambda\right)  ,$ are real coefficients polynomials for
$\lambda.$

\section{\textbf{Operator Formulation and Properties of Spectrum}}

In this section, we present the space\newline$H:=L_{2}(\Lambda)\oplus%
\mathbb{C}
^{m_{a}}\oplus%
\mathbb{C}
^{m_{b}}\oplus%
{\displaystyle\sum\limits_{i=1}^{n}}
\mathbb{C}
^{r_{i}}$ where $m_{a}=\max\left\{  \deg a_{i}(\lambda)\right\}  ,$
$m_{b}=\max\left\{  \deg b_{i}(\lambda)\right\}  ,$ $r_{i}=\max\left\{
\deg\gamma_{i}(\lambda)\right\}  .$ We define the norm on space $H$ by%
\begin{align}
\left\Vert Y\right\Vert ^{2}  &  :=%
{\displaystyle\int\limits_{a}^{b}}
\rho\left(  x\right)  \left\{  \left\vert y_{1}(x)\right\vert ^{2}+\left\vert
y_{2}(x)\right\vert ^{2}\right\}  dx+%
{\displaystyle\sum\limits_{i=1}^{m_{a}}}
Y_{i}^{1}\overline{Y_{i}^{1}}\tag{6}\\
&  \text{ \ \ \ \ \ \ \ \ \ \ \ \ }+%
{\displaystyle\sum\limits_{j=1}^{m_{b}}}
Y_{j}^{2}\overline{Y_{j}^{2}}+%
{\displaystyle\sum\limits_{j=1}^{n}}
{\displaystyle\sum\limits_{i=1}^{n}}
Y_{i}^{3_{j}}\overline{Y_{i}^{3_{j}}}\nonumber
\end{align}

for

$Y\in H,$ $Y=\left(  f(x),Y^{1},Y^{2},Y^{3}\right)  $ $,$ $Y^{1}=\left(
Y_{1}^{1},Y_{2}^{1},\ldots,Y_{m_{a}}^{1}\right)  ,$ $Y^{2}=\left(  Y_{1}%
^{2},Y_{2}^{2},\ldots,Y_{m_{b}}^{2}\right)  ,$ $Y^{3}=\left(  Y_{1}^{3_{i}%
},Y_{2}^{3_{i}},\ldots,Y_{r_{i}}^{3_{i}}\right)  $ $\left(  i=\overline
{1,n}\right)  .\smallskip$

Consider the operator $T$ defined by the domain$\smallskip$\newline%
$D(T)=\left\{  Y\in H:%
\begin{array}
[c]{l}%
i)\text{ }y(x)=\left(
\begin{array}
[c]{c}%
y_{1}(x)\\
y_{2}(x)
\end{array}
\right)  \in AC(\Lambda),\text{ }lf\in L_{2}(\Lambda)\medskip\\
ii)\text{ }y_{1}(\xi_{i}+0)-\theta_{i}y_{1}(\xi_{i}-0)=0\medskip\\
iii)\text{ }Y_{1}^{1}=a_{m_{a}2}y_{2}(a)-a_{m_{a}1}y_{1}(a)\medskip\\
iv)\text{ }Y_{1}^{2}=b_{m_{b}2}y_{2}(b)-b_{m_{b}1}y_{1}(b)\medskip\\
v)\text{ }Y_{1}^{3_{i}}=\gamma_{r_{i}i}y_{1}(\xi_{i}-0)\medskip
\end{array}
\right\}  $\smallskip$\medskip$\newline such that$\smallskip$\newline%
$TY:=W=(ly,W^{1},W^{2},W^{3_{\dot{I}}}),$ $W^{1}=\left(  W_{1}^{1},W_{2}%
^{1},\ldots,W_{m_{a}}^{1}\right)  ,\medskip$

$W^{2}=\left(  W_{1}^{2},W_{2}^{2},\ldots,W_{m_{b}}^{2}\right)  ,$
$W^{3_{\dot{I}}}=\left(  W_{1}^{3_{i}},W_{2}^{3_{i}},\ldots,W_{r_{i}}^{3_{i}%
},W_{r_{n+1}}^{3}\right)  $ $\left(  i=\overline{1,n}\right)  \medskip$

$W_{i}^{1}=Y_{i+1}^{1}-a_{m_{a}-i,2}y_{2}(a)+a_{m_{a}-i,1}y_{1}(a),$
$i=\overline{1,m_{a}-1}\medskip$

$W_{m_{a}}^{1}=-a_{02}y_{2}(a)+a_{01}y_{1}(a)\medskip$

$W_{j}^{2}=Y_{j+1}^{1}-b_{m_{b}-j,2}y_{2}(b)+b_{m_{b}-j,1}y_{1}(b),$
$j=\overline{1,m_{b}-1}\medskip$

$W_{m_{b}}^{2}=-b_{02}y_{2}(b)+b_{01}y_{1}(b)\medskip$

$W_{1}^{3}=Y_{K+1}^{3_{1}}-\gamma_{r_{1}-k,1}y_{1}(\xi_{1}-0),$ $k=\overline
{1,r_{1}-1}\medskip$

$W_{2}^{3}=Y_{K+1}^{3_{2}}-\gamma_{r_{2}-k,2}y_{1}(\xi_{2}-0),$ $k=\overline
{1,r_{2}-1}\medskip$

$\vdots\medskip$

$W_{r_{n}}^{3}=Y_{K+1}^{3_{n}}-\gamma_{r_{n}-k,n}y_{1}(\xi_{n}-0),$
$k=\overline{1,r_{n}-1}\medskip$

$W_{r_{n+1}}^{3}=-\gamma_{0t}y_{1}(\xi_{t}-0)+y_{2}(\xi_{t}+0)-\theta_{i}%
^{-1}y_{2}(\xi_{t}-0),\left(  t=\overline{1,n}\right)  \medskip$

Thus, we can rewrite the considered problem (1)-(5) in the operator form as

$TY=\lambda Y$. \smallskip

We define the solutions\medskip

$\varphi(x,\lambda)=\varphi_{i}(x,\lambda),$ $x\in\left(  \xi_{i},\xi
_{i+1}\right)  ;\ \psi(x,\lambda)=\psi_{i}(x,\lambda),$ $x\in\left(  \xi
_{i},\xi_{i+1}\right)  $\smallskip\ $\left(  i=\overline{0,n}\right)
\medskip$

$\varphi_{1}(x,\lambda)=\left(  \varphi_{i1}(x,\lambda),\varphi_{i2}%
(x,\lambda)\right)  ^{T}$and $\psi_{1}(x,\lambda)=\left(  \psi_{i1}%
(x,\lambda),\psi_{i2}(x,\lambda)\right)  ^{T},$ $\left(  i=\overline
{0,n}\right)  $ of equation (1) by the initial conditions%
\begin{equation}
\varphi_{01}(a,\lambda)=a_{2}\left(  \lambda\right)  ,\ \varphi_{02}%
(a,\lambda)=a_{1}\left(  \lambda\right)  \tag{7}%
\end{equation}

and similarly;%
\begin{equation}
\psi_{n1}(b,\lambda)=b_{2}\left(  \lambda\right)  ,\ \psi_{n2}(b,\lambda
)=b_{1}\left(  \lambda\right)  \tag{8}%
\end{equation}

respectively.

\textbf{Lemma 1 }The following asymptotic behaviours hold for $\left\vert
\lambda\right\vert \rightarrow\infty$:$\smallskip$\newline$\qquad\varphi
_{01}(x,\lambda)=\left\{
\begin{array}
[c]{l}%
a_{m_{2}2}\lambda^{m_{2}}\cos\lambda\rho_{0}(x-a)+o(\lambda^{m_{2}}%
\exp\left\vert \operatorname{Im}\lambda\right\vert \left(  x-a\right)
\rho_{0},\text{ }\deg a_{2}\left(  \lambda\right)  >\deg a_{1}\left(
\lambda\right)  \medskip\\
-a_{m_{1}1}\lambda^{m_{1}}\sin\lambda\rho_{0}(x-a)+o(\lambda^{m_{1}}%
\exp\left\vert \operatorname{Im}\lambda\right\vert \left(  x-a\right)
\rho_{0},\text{ }\deg a_{1}\left(  \lambda\right)  >\deg a_{2}\left(
\lambda\right)  \medskip\\
\lambda^{m}(a_{m_{2}2}\cos\lambda\rho_{0}(x-a)-a_{m_{1}1}\sin\lambda\rho
_{0}(x-a))\medskip\\
+o(\lambda^{m}\exp\left\vert \operatorname{Im}\lambda\right\vert \left(
x-a\right)  \rho_{0},\text{ }\deg a_{1}\left(  \lambda\right)  =\deg
a_{2}\left(  \lambda\right)  \medskip
\end{array}
\right.  \smallskip$

$\varphi_{02}(x,\lambda)=\left\{
\begin{array}
[c]{l}%
a_{m_{2}2}\lambda^{m_{2}}\sin\lambda\rho_{0}(x-a)+o(\lambda^{m_{2}}%
\exp\left\vert \operatorname{Im}\lambda\right\vert \left(  x-a\right)
\rho_{0},\text{ }\deg a_{2}\left(  \lambda\right)  >\deg a_{1}\left(
\lambda\right)  \medskip\\
a_{m_{1}1}\lambda^{m_{1}}\cos\lambda\rho_{0}(x-a)+o(\lambda^{m_{1}}%
\exp\left\vert \operatorname{Im}\lambda\right\vert \left(  x-a\right)
\rho_{0},\text{ }\deg a_{1}\left(  \lambda\right)  >\deg a_{2}\left(
\lambda\right)  \medskip\\
\lambda^{m}(a_{m_{1}1}\cos\lambda\rho_{0}(x-a)+a_{m_{2}2}\sin\lambda\rho
_{0}(x-a))\medskip\\
+o(\lambda^{m}\exp\left\vert \operatorname{Im}\lambda\right\vert \left(
x-a\right)  \rho_{0},\text{ }\deg a_{1}\left(  \lambda\right)  =\deg
a_{2}\left(  \lambda\right)  \medskip
\end{array}
\right.  \smallskip$

$\varphi_{11}(x,\lambda)=\left\{
\begin{array}
[c]{l}%
\begin{array}
[c]{c}%
-\gamma_{r_{1}1}a_{m_{2}2}\lambda^{r_{1}+m_{2}}\cos\lambda\rho_{0}(\xi
_{1}-a)\sin\lambda\rho_{1}\left(  x-\xi_{1}\right)  \medskip\\
+o(\lambda^{r_{1}+m_{2}}\exp\left\vert \operatorname{Im}\lambda\right\vert
(\left(  \xi_{1}-a\right)  \rho_{0}+\left(  x-\xi_{1}\right)  \rho_{1}),\text{
}\deg a_{2}\left(  \lambda\right)  >\deg a_{1}\left(  \lambda\right)  \medskip
\end{array}
\smallskip\\%
\begin{array}
[c]{c}%
\gamma_{r_{1}1}a_{m_{1}1}\lambda^{r_{1}+m_{1}}\sin\lambda\rho_{0}(\xi
_{1}-a)\sin\lambda\rho_{1}\left(  x-\xi_{1}\right)  \medskip\\
+o(\lambda^{r_{1}+m_{1}}\exp\left\vert \operatorname{Im}\lambda\right\vert
(\left(  \xi_{1}-a\right)  \rho_{0}+\left(  x-\xi_{1}\right)  \rho_{1}),\text{
}\deg a_{1}\left(  \lambda\right)  >\deg a_{2}\left(  \lambda\right)  \medskip
\end{array}
\smallskip\\
\multicolumn{1}{c}{%
\begin{array}
[c]{l}%
\gamma_{r_{1}1}\lambda^{r_{1}+m}(-a_{m_{2}2}\cos\lambda\rho_{0}(\xi_{1}%
-a)\sin\lambda\rho_{1}\left(  x-\xi_{1}\right)  \medskip\\
+a_{m_{1}1}\sin\lambda\rho_{0}(\xi_{1}-a)\sin\lambda\rho_{1}\left(  x-\xi
_{1}\right)  \medskip\\
+o(\lambda^{r_{1}+m}\exp\left\vert \operatorname{Im}\lambda\right\vert
(\left(  \xi_{1}-a\right)  \rho_{0}+\left(  x-\xi_{1}\right)  \rho_{1}),\text{
}\deg a_{1}\left(  \lambda\right)  =\deg a_{2}\left(  \lambda\right)  \medskip
\end{array}
\smallskip}%
\end{array}
\right.  \smallskip$

$\varphi_{11}(x,\lambda)=\left\{
\begin{array}
[c]{l}%
\begin{array}
[c]{l}%
-\gamma_{r_{1}1}a_{m_{2}2}\lambda^{r_{1}+m_{2}}\cos\lambda\rho_{0}(\xi
_{1}-a)\sin\lambda\rho_{1}\left(  x-\xi_{1}\right)  \medskip\\
+o(\lambda^{r_{1}+m_{2}}\exp\left\vert \operatorname{Im}\lambda\right\vert
(\left(  \xi_{1}-a\right)  \rho_{0}+\left(  x-\xi_{1}\right)  \rho_{1}),\deg
a_{2}\left(  \lambda\right)  >\deg a_{1}\left(  \lambda\right)  \medskip
\end{array}
\smallskip\\%
\begin{array}
[c]{l}%
\gamma_{r_{1}1}a_{m_{1}1}\lambda^{r_{1}+m_{1}}\sin\lambda\rho_{0}(\xi
_{1}-a)\sin\lambda\rho_{1}\left(  x-\xi_{1}\right)  \medskip\\
+o(\lambda^{r_{1}+m_{1}}\exp\left\vert \operatorname{Im}\lambda\right\vert
(\left(  \xi_{1}-a\right)  \rho_{0}+\left(  x-\xi_{1}\right)  \rho_{1}),\deg
a_{1}\left(  \lambda\right)  >\deg a_{2}\left(  \lambda\right)  \medskip
\end{array}
\smallskip\\%
\begin{array}
[c]{l}%
\gamma_{r_{1}1}\lambda^{r_{1}+m}(-a_{m_{2}2}\cos\lambda\rho_{0}(\xi_{1}%
-a)\sin\lambda\rho_{1}\left(  x-\xi_{1}\right)  +a_{m_{1}1}\sin\lambda\rho
_{0}(\xi_{1}-a)\sin\lambda\rho_{1}\left(  x-\xi_{1}\right)  \medskip\\
+o(\lambda^{r_{1}+m}\exp\left\vert \operatorname{Im}\lambda\right\vert
(\left(  \xi_{1}-a\right)  \rho_{0}+\left(  x-\xi_{1}\right)  \rho_{1}),\deg
a_{1}\left(  \lambda\right)  =\deg a_{2}\left(  \lambda\right)  \medskip
\end{array}
\smallskip
\end{array}
\right.  \smallskip$

$\varphi_{12}(x,\lambda)=\left\{
\begin{array}
[c]{l}%
\begin{array}
[c]{l}%
\gamma_{r_{1}1}a_{m_{2}2}\lambda^{r_{1}+m_{2}}\cos\lambda\rho_{0}(\xi
_{1}-a)\cos\lambda\rho_{1}\left(  x-\xi_{1}\right)  \medskip\\
+o(\lambda^{r_{1}+m_{2}}\exp\left\vert \operatorname{Im}\lambda\right\vert
(\left(  \xi_{1}-a\right)  \rho_{0}+\left(  x-\xi_{1}\right)  \rho_{1}),\text{
}\deg a_{2}\left(  \lambda\right)  >\deg a_{1}\left(  \lambda\right)  \medskip
\end{array}
\smallskip\\%
\begin{array}
[c]{l}%
-\gamma_{r_{1}1}a_{m_{1}1}\lambda^{r_{1}+m_{1}}\sin\lambda\rho_{0}(\xi
_{1}-a)\cos\lambda\rho_{1}\left(  x-\xi_{1}\right)  \medskip\\
+o(\lambda^{r_{1}+m_{1}}\exp\left\vert \operatorname{Im}\lambda\right\vert
(\left(  \xi_{1}-a\right)  \rho_{0}+\left(  x-\xi_{1}\right)  \rho_{1}),\text{
}\deg a_{1}\left(  \lambda\right)  >\deg a_{2}\left(  \lambda\right)  \medskip
\end{array}
\smallskip\\%
\begin{array}
[c]{l}%
\gamma_{r_{1}1}\lambda^{r_{1}+m}(a_{m_{2}2}\cos\lambda\rho_{0}(\xi_{1}%
-a)\cos\lambda\rho_{1}\left(  x-\xi_{1}\right)  a_{m_{1}1}\sin\lambda\rho
_{0}(\xi_{1}-a)\cos\lambda\rho_{1}\left(  x-\xi_{1}\right)  \medskip\\
+o(\lambda^{r_{1}+m}\exp\left\vert \operatorname{Im}\lambda\right\vert
(\left(  \xi_{1}-a\right)  \rho_{0}+\left(  x-\xi_{1}\right)  \rho_{1}),\text{
}\deg a_{1}\left(  \lambda\right)  =\deg a_{2}\left(  \lambda\right)  \medskip
\end{array}
\smallskip
\end{array}
\right.  \smallskip$

$\varphi_{21}(x,\lambda)=\left\{
\begin{array}
[c]{l}%
\begin{array}
[c]{l}%
\gamma_{r_{1}1}\gamma_{r_{2}2}a_{m_{2}2}\lambda^{r_{1}+r_{2}+m_{2}}\cos
\lambda\rho_{0}(\xi_{1}-a)\sin\lambda\rho_{1}\left(  \xi_{2}-\xi_{1}\right)
\sin\lambda\rho_{2}\left(  x-\xi_{2}\right)  \medskip\\
+o(\lambda^{r_{1}+r_{2}+m_{2}}\exp\left\vert \operatorname{Im}\lambda
\right\vert (\left(  \xi_{1}-a\right)  \rho_{0}+\left(  \xi_{2}-\xi
_{1}\right)  \rho_{1}+\left(  x-\xi_{2}\right)  \rho_{2}),\text{ }\deg
a_{2}\left(  \lambda\right)  >\deg a_{1}\left(  \lambda\right)  \medskip
\end{array}
\smallskip\\%
\begin{array}
[c]{l}%
-\gamma_{r_{1}1}\gamma_{r_{2}2}a_{m_{1}1}\lambda^{r_{1}+r_{2}+m_{1}}%
\sin\lambda\rho_{0}(\xi_{1}-a)\sin\lambda\rho_{1}\left(  \xi_{2}-\xi
_{1}\right)  \sin\lambda\rho_{2}\left(  x-\xi_{2}\right)  \medskip\\
+o(\lambda^{r_{1}+r_{2}+m_{1}}\exp\left\vert \operatorname{Im}\lambda
\right\vert (\left(  \xi_{1}-a\right)  \rho_{0}+\left(  \xi_{2}-\xi
_{1}\right)  \rho_{1}+\left(  x-\xi_{2}\right)  \rho_{2}),\text{ }\deg
a_{1}\left(  \lambda\right)  >\deg a_{2}\left(  \lambda\right)  \medskip
\end{array}
\smallskip\\%
\begin{array}
[c]{l}%
\gamma_{r_{1}1}\gamma_{r_{2}2}\lambda^{r_{1}+r_{2}+m}(a_{m_{2}2}\cos
\lambda\rho_{0}(\xi_{1}-a)\sin\lambda\rho_{1}\left(  \xi_{2}-\xi_{1}\right)
\sin\lambda\rho_{2}\left(  x-\xi_{2}\right)  \medskip\\
-a_{m_{1}1}\sin\lambda\rho_{0}(\xi_{1}-a)\sin\lambda\rho_{1}\left(  \xi
_{2}-\xi_{1}\right)  \sin\lambda\rho_{2}\left(  x-\xi_{2}\right)  )\medskip\\
+o(\lambda^{r_{1}+r_{2}+m}\exp\left\vert \operatorname{Im}\lambda\right\vert
(\left(  \xi_{1}-a\right)  \rho_{0}+\left(  \xi_{2}-\xi_{1}\right)  \rho
_{1}+\left(  x-\xi_{2}\right)  \rho_{2}),\text{ }\deg a_{1}\left(
\lambda\right)  =\deg a_{2}\left(  \lambda\right)  \medskip
\end{array}
\smallskip
\end{array}
\right.  \smallskip$

$\varphi_{22}(x,\lambda)=\left\{
\begin{array}
[c]{l}%
\begin{array}
[c]{l}%
-\gamma_{r_{1}1}\gamma_{r_{2}2}a_{m_{2}2}\lambda^{r_{1}+r_{2}+m_{2}}%
\cos\lambda\rho_{0}(\xi_{1}-a)\sin\lambda\rho_{1}\left(  \xi_{2}-\xi
_{1}\right)  \cos\lambda\rho_{2}\left(  x-\xi_{2}\right)  \medskip\\
+o(\lambda^{r_{1}+r_{2}+m_{2}}\exp\left\vert \operatorname{Im}\lambda
\right\vert (\left(  \xi_{1}-a\right)  \rho_{0}+\left(  \xi_{2}-\xi
_{1}\right)  \rho_{1}+\left(  x-\xi_{2}\right)  \rho_{2}),\text{ }\deg
a_{2}\left(  \lambda\right)  >\deg a_{1}\left(  \lambda\right)  \medskip
\end{array}
\smallskip\\%
\begin{array}
[c]{l}%
\gamma_{r_{1}1}\gamma_{r_{2}2}a_{m_{1}1}\lambda^{r_{1}+r_{2}+m_{1}}\sin
\lambda\rho_{0}(\xi_{1}-a)\sin\lambda\rho_{1}\left(  \xi_{2}-\xi_{1}\right)
\cos\lambda\rho_{2}\left(  x-\xi_{2}\right)  \medskip\\
+o(\lambda^{r_{1}+r_{2}+m_{1}}\exp\left\vert \operatorname{Im}\lambda
\right\vert (\left(  \xi_{1}-a\right)  \rho_{0}+\left(  \xi_{2}-\xi
_{1}\right)  \rho_{1}+\left(  x-\xi_{2}\right)  \rho_{2}),\text{ }\deg
a_{1}\left(  \lambda\right)  >\deg a_{2}\left(  \lambda\right)  \medskip
\end{array}
\smallskip\\%
\begin{array}
[c]{l}%
\gamma_{r_{1}1}\gamma_{r_{2}2}\lambda^{r_{1}+r_{2}+m}(-a_{m_{2}2}\cos
\lambda\rho_{0}(\xi_{1}-a)\sin\lambda\rho_{1}\left(  \xi_{2}-\xi_{1}\right)
\cos\lambda\rho_{2}\left(  x-\xi_{2}\right)  \medskip\\
+a_{m_{1}1}\sin\lambda\rho_{0}(\xi_{1}-a)\sin\lambda\rho_{1}\left(  \xi
_{2}-\xi_{1}\right)  \cos\lambda\rho_{2}\left(  x-\xi_{2}\right)  \medskip\\
+o(\lambda^{r_{1}+r_{2}+m}\exp\left\vert \operatorname{Im}\lambda\right\vert
(\left(  \xi_{1}-a\right)  \rho_{0}+\left(  \xi_{2}-\xi_{1}\right)  \rho
_{1}+\left(  x-\xi_{2}\right)  \rho_{2}),\text{ }\deg a_{1}\left(
\lambda\right)  =\deg a_{2}\left(  \lambda\right)  \medskip
\end{array}
\smallskip
\end{array}
\right.  \medskip$

$,\ldots,$

$\varphi_{n1}(x,\lambda)=\left\{
\begin{array}
[c]{l}%
\begin{array}
[c]{l}%
\left(  -1\right)  ^{n}a_{m_{2}2}\lambda^{A+m_{2}}\times\left(
{\displaystyle\prod\limits_{i=1}^{n}}
\gamma_{r_{i}i}\right)  \times\cos\lambda\rho_{0}(\xi_{1}-a)\sin\lambda
\rho_{n}\left(  x-\xi_{n}\right)  \times\medskip\\
\times\left(
{\displaystyle\prod\limits_{i=2}^{n}}
\sin\lambda\rho_{i-1}\left(  \xi_{i}-\xi_{i-1}\right)  \right)  \medskip\\
+o(\lambda^{A+m_{2}}\exp\left\vert \operatorname{Im}\lambda\right\vert
(\left(  b-\xi_{n}\right)  \rho_{n}+%
{\displaystyle\sum\limits_{i=1}^{n}}
\left(  \xi_{i}-\xi_{i-1}\right)  \rho_{i-1}),\text{ }\deg a_{2}\left(
\lambda\right)  >\deg a_{1}\left(  \lambda\right)  \medskip
\end{array}
\smallskip\\%
\begin{array}
[c]{l}%
\left(  -1\right)  ^{n+1}a_{m_{1}1}\lambda^{A+m_{1}}\times\left(
{\displaystyle\prod\limits_{i=1}^{n}}
\gamma_{r_{i}i}\right)  \sin\lambda\rho_{0}(\xi_{1}-a)\sin\lambda\rho
_{n}\left(  x-\xi_{n}\right)  \times\\
\times\left(
{\displaystyle\prod\limits_{i=2}^{n}}
\sin\lambda\rho_{i-1}\left(  \xi_{i}-\xi_{i-1}\right)  \right)  \medskip\\
+o(\lambda^{A+m_{1}}\exp\left\vert \operatorname{Im}\lambda\right\vert
(\left(  b-\xi_{n}\right)  \rho_{n}+%
{\displaystyle\sum\limits_{i=1}^{n}}
\left(  \xi_{i}-\xi_{i-1}\right)  \rho_{i-1}),\text{ }\deg a_{1}\left(
\lambda\right)  >\deg a_{2}\left(  \lambda\right)  \medskip
\end{array}
\smallskip\\
\lambda^{A+m}\times\left(
{\displaystyle\prod\limits_{i=1}^{n}}
\gamma_{r_{i}i}\right)  \times\\
\times\left\{
\begin{array}
[c]{c}%
\left(  -1\right)  ^{n}a_{m_{2}2}\cos\lambda\rho_{0}(\xi_{1}-a)\sin\lambda
\rho_{n}\left(  x-\xi_{n}\right)  \times\left(
{\displaystyle\prod\limits_{i=2}^{n}}
\sin\lambda\rho_{i-1}\left(  \xi_{i}-\xi_{i-1}\right)  \right)  \medskip\\
+\left(  -1\right)  ^{n+1}a_{m_{1}1}\sin\lambda\rho_{0}(\xi_{1}-a)\sin
\lambda\rho_{n}\left(  x-\xi_{n}\right)  \times\left(
{\displaystyle\prod\limits_{i=2}^{n}}
\sin\lambda\rho_{i-1}\left(  \xi_{i}-\xi_{i-1}\right)  \right)  \medskip
\end{array}
\right\}  \smallskip\\
+o(\lambda^{A+m}\exp\left\vert \operatorname{Im}\lambda\right\vert (\left(
b-\xi_{n}\right)  \rho_{n}+%
{\displaystyle\sum\limits_{i=1}^{n}}
\left(  \xi_{i}-\xi_{i-1}\right)  \rho_{i-1}),\text{ }\deg a_{1}\left(
\lambda\right)  =\deg a_{2}\left(  \lambda\right)
\end{array}
\right.  \smallskip$

$\varphi_{n2}(x,\lambda)=\left\{
\begin{array}
[c]{l}%
\begin{array}
[c]{l}%
\left(  -1\right)  ^{n+1}a_{m_{2}2}\lambda^{A+m_{2}}\times\left(
{\displaystyle\prod\limits_{i=1}^{n}}
\gamma_{r_{i}i}\right)  \cos\lambda\rho_{0}(\xi_{1}-a)\cos\lambda\rho
_{n}\left(  x-\xi_{n}\right)  \times\\
\times\left(
{\displaystyle\prod\limits_{i=2}^{n}}
\sin\lambda\rho_{i-1}\left(  \xi_{i}-\xi_{i-1}\right)  \right)  \medskip\\
+o(\lambda^{A+m_{2}}\exp\left\vert \operatorname{Im}\lambda\right\vert
(\left(  b-\xi_{n}\right)  \rho_{n}+%
{\displaystyle\sum\limits_{i=1}^{n}}
\left(  \xi_{i}-\xi_{i-1}\right)  \rho_{i-1}),\text{ }\deg a_{2}\left(
\lambda\right)  >\deg a_{1}\left(  \lambda\right)  \medskip
\end{array}
\smallskip\\%
\begin{array}
[c]{l}%
\left(  -1\right)  ^{n}a_{m_{1}1}\lambda^{A+m_{1}}\times\left(
{\displaystyle\prod\limits_{i=1}^{n}}
\gamma_{r_{i}i}\right)  \sin\lambda\rho_{0}(\xi_{1}-a)\cos\lambda\rho
_{n}\left(  x-\xi_{n}\right)  \times\\
\times\left(
{\displaystyle\prod\limits_{i=2}^{n}}
\sin\lambda\rho_{i-1}\left(  \xi_{i}-\xi_{i-1}\right)  \right)  \medskip\\
+o(\lambda^{A+m_{1}}\exp\left\vert \operatorname{Im}\lambda\right\vert
(\left(  b-\xi_{n}\right)  \rho_{n}+%
{\displaystyle\sum\limits_{i=1}^{n}}
\left(  \xi_{i}-\xi_{i-1}\right)  \rho_{i-1}),\text{ }\deg a_{1}\left(
\lambda\right)  >\deg a_{2}\left(  \lambda\right)  \medskip
\end{array}
\smallskip\\%
\begin{array}
[c]{l}%
\lambda^{A+m}\times\left(
{\displaystyle\prod\limits_{i=1}^{n}}
\gamma_{r_{i}i}\right)  \times\\
\times\left\{
\begin{array}
[c]{c}%
\left(  -1\right)  ^{n+1}a_{m_{2}2}\cos\lambda\rho_{0}(\xi_{1}-a)\cos
\lambda\rho_{n}\left(  x-\xi_{n}\right)  \times\left(
{\displaystyle\prod\limits_{i=2}^{n}}
\sin\lambda\rho_{i-1}\left(  \xi_{i}-\xi_{i-1}\right)  \right)  \medskip\\
+\left(  -1\right)  ^{n}a_{m_{1}1}\sin\lambda\rho_{0}(\xi_{1}-a)\cos
\lambda\rho_{n}\left(  x-\xi_{n}\right)  \times\left(
{\displaystyle\prod\limits_{i=2}^{n}}
\sin\lambda\rho_{i-1}\left(  \xi_{i}-\xi_{i-1}\right)  \right)  \medskip
\end{array}
\right\}  \medskip\\
+o(\lambda^{A+m}\exp\left\vert \operatorname{Im}\lambda\right\vert (\left(
b-\xi_{n}\right)  \rho_{n}+%
{\displaystyle\sum\limits_{i=1}^{n}}
\left(  \xi_{i}-\xi_{i-1}\right)  \rho_{i-1}),\text{ }\deg a_{1}\left(
\lambda\right)  =\deg a_{2}\left(  \lambda\right)  \medskip
\end{array}
\smallskip
\end{array}
\right.  \smallskip$

where $A=r_{1}+r_{2}+\cdots+r_{n}.$

\textbf{Lemma 2 }The following asymptotic behaviours hold for $\left\vert
\lambda\right\vert \rightarrow\infty$:$\smallskip$

$\psi_{n1}(x,\lambda)=\left\{
\begin{array}
[c]{l}%
b_{m_{4}2}\lambda^{m_{4}}\cos\lambda\rho_{n}(x-b)+o(\lambda^{m_{4}}%
\exp\left\vert \operatorname{Im}\lambda\right\vert \left(  x-b\right)
\rho_{n},\text{ }\deg b_{2}\left(  \lambda\right)  >\deg b_{1}\left(
\lambda\right)  \medskip\\
-b_{m_{3}1}\lambda^{m_{3}}\sin\lambda\rho_{n}(x-b)+o(\lambda^{m_{3}}%
\exp\left\vert \operatorname{Im}\lambda\right\vert \left(  x-b\right)
\rho_{n},\text{ }\deg b_{1}\left(  \lambda\right)  >\deg b_{2}\left(
\lambda\right)  \medskip\\
\lambda^{m_{b}}(b_{m_{4}2}\cos\lambda\rho_{n}(x-b)-b_{m_{3}1}\sin\lambda
\rho_{n}(x-b))\medskip\\
+o(\lambda^{m_{b}}\exp\left\vert \operatorname{Im}\lambda\right\vert \left(
x-b\right)  \rho_{n}),\text{ }\deg b_{1}\left(  \lambda\right)  =\deg
b_{2}\left(  \lambda\right)  \medskip
\end{array}
\right.  \smallskip$

$\psi_{n2}(x,\lambda)=\left\{
\begin{array}
[c]{l}%
b_{m_{4}2}\lambda^{m_{4}}\sin\lambda\rho_{n}(x-b)+o(\lambda^{m_{4}}%
\exp\left\vert \operatorname{Im}\lambda\right\vert \left(  x-b\right)
\rho_{n},\text{ }\deg b_{2}\left(  \lambda\right)  >\deg b_{1}\left(
\lambda\right)  \medskip\\
b_{m_{3}1}\lambda^{m_{3}}\cos\lambda\rho_{n}(x-b)+o(\lambda^{m_{3}}%
\exp\left\vert \operatorname{Im}\lambda\right\vert \left(  x-b\right)
\rho_{n},\text{ }\deg b_{1}\left(  \lambda\right)  >\deg b_{2}\left(
\lambda\right)  \medskip\\
\lambda^{m_{b}}(b_{m_{4}2}\sin\lambda\rho_{n}(x-b)+b_{m_{3}1}\cos\lambda
\rho_{n}(x-b))\medskip\\
+o(\lambda^{m_{b}}\exp\left\vert \operatorname{Im}\lambda\right\vert \left(
x-b\right)  \rho_{n}),\text{ }\deg b_{1}\left(  \lambda\right)  =\deg
b_{2}\left(  \lambda\right)  \medskip
\end{array}
\right.  \smallskip$

$\psi_{n-1,1}(x,\lambda)=\left\{
\begin{array}
[c]{c}%
\begin{array}
[c]{l}%
b_{m_{4}2}\gamma_{r_{n}n}\lambda^{r_{n}+m_{4}}\cos\lambda\rho_{n}(\xi
_{n}-b)\sin\lambda\rho_{n-1}\left(  x-\xi_{n}\right)  \medskip\\
+o(\lambda^{r_{n}+m_{4}}\exp\left\vert \operatorname{Im}\lambda\right\vert
(\left(  \xi_{n}-b\right)  \rho_{n}+\left(  x-\xi_{n}\right)  \rho_{n-1}),\deg
b_{2}\left(  \lambda\right)  >\deg b_{1}\left(  \lambda\right)  \medskip
\end{array}
\\%
\begin{array}
[c]{l}%
-b_{m_{3}1}\gamma_{r_{n}n}\lambda^{r_{n}+m_{3}}\sin\lambda\rho_{n}(\xi
_{n}-b)\sin\lambda\rho_{n-1}\left(  x-\xi_{n}\right)  \medskip\\
+o(\lambda^{r_{n}+m_{3}}\exp\left\vert \operatorname{Im}\lambda\right\vert
(\left(  \xi_{n}-b\right)  \rho_{n}+\left(  x-\xi_{n}\right)  \rho
_{n-1}),\text{ }\deg b_{1}\left(  \lambda\right)  >\deg b_{2}\left(
\lambda\right)  \medskip
\end{array}
\\%
\begin{array}
[c]{l}%
\gamma_{r_{n}n}\lambda^{r_{n}+m_{b}}\left(  -b_{m_{3}1}\sin\lambda\rho_{n}%
(\xi_{n}-b)+b_{m_{4}2}\cos\lambda\rho_{n}(\xi_{n}-b)\right)  \sin\lambda
\rho_{n-1}\left(  x-\xi_{n}\right)  \medskip\\
+o(\lambda^{r_{n}+m_{b}}\exp\left\vert \operatorname{Im}\lambda\right\vert
(\left(  \xi_{n}-b\right)  \rho_{n}+\left(  x-\xi_{n}\right)  \rho
_{n-1}),\text{ }\deg b_{1}\left(  \lambda\right)  =\deg b_{2}\left(
\lambda\right)  \medskip
\end{array}
\end{array}
\right.  \smallskip$

$\psi_{n-1,2}(x,\lambda)=\left\{
\begin{array}
[c]{l}%
\begin{array}
[c]{l}%
-b_{m_{4}2}\gamma_{r_{n}n}\lambda^{r_{n}+m_{4}}\cos\lambda\rho_{n}(\xi
_{n}-b)\cos\lambda\rho_{n-1}\left(  x-\xi_{n}\right)  \medskip\\
+o(\lambda^{r_{n}+m_{4}}\exp\left\vert \operatorname{Im}\lambda\right\vert
(\left(  \xi_{n}-b\right)  \rho_{n}+\left(  x-\xi_{n}\right)  \rho
_{n-1}),\text{ }\deg b_{2}\left(  \lambda\right)  >\deg b_{1}\left(
\lambda\right)  \medskip
\end{array}
\\%
\begin{array}
[c]{l}%
b_{m_{3}1}\gamma_{r_{n}n}\lambda^{r_{n}+m_{3}}\sin\lambda\rho_{n}(\xi
_{n}-b)\cos\lambda\rho_{n-1}\left(  x-\xi_{n}\right)  \medskip\\
+o(\lambda^{r_{n}+m_{3}}\exp\left\vert \operatorname{Im}\lambda\right\vert
(\left(  \xi_{n}-b\right)  \rho_{n}+\left(  x-\xi_{n}\right)  \rho
_{n-1}),\text{ }\deg b_{1}\left(  \lambda\right)  >\deg b_{2}\left(
\lambda\right)  \medskip
\end{array}
\\%
\begin{array}
[c]{l}%
\gamma_{r_{n}n}\lambda^{r_{n}+m_{b}}\left(  b_{m_{3}1}\sin\lambda\rho_{n}%
(\xi_{n}-b)-b_{m_{4}2}\cos\lambda\rho_{n}(\xi_{n}-b)\right)  \cos\lambda
\rho_{n-1}\left(  x-\xi_{n}\right)  \medskip\\
+o(\lambda^{r_{n}+m_{b}}\exp\left\vert \operatorname{Im}\lambda\right\vert
(\left(  \xi_{n}-b\right)  \rho_{n}+\left(  x-\xi_{n}\right)  \rho
_{n-1}),\text{ }\deg b_{1}\left(  \lambda\right)  =\deg b_{2}\left(
\lambda\right)  \medskip
\end{array}
\end{array}
\right.  \smallskip$

$\psi_{n-2,1}(x,\lambda)=\left\{
\begin{array}
[c]{l}%
\begin{array}
[c]{l}%
-b_{m_{4}2}\gamma_{r_{n}n}\gamma_{r_{n-1}n-1}\lambda^{r_{n}+r_{n-1}+m_{4}%
}\times\medskip\\
\times\cos\lambda\rho_{n}(\xi_{n}-b)\sin\lambda\rho_{n-1}\left(  \xi_{n}%
-\xi_{n-1}\right)  \sin\lambda\rho_{n-2}\left(  x-\xi_{n-1}\right)  \medskip\\
+o(\lambda^{r_{n}+r_{n-1}+m_{4}}\exp\left\vert \operatorname{Im}%
\lambda\right\vert ((\xi_{n}-b)\rho_{n}+\left(  \xi_{n}-\xi_{n-1}\right)
\rho_{n-1}+\left(  x-\xi_{n-1}\right)  \rho_{n-2})\\
,\text{ }\deg b_{2}\left(  \lambda\right)  >\deg b_{1}\left(  \lambda\right)
\medskip
\end{array}
\\%
\begin{array}
[c]{l}%
b_{m_{3}1}\gamma_{r_{n}n}\gamma_{r_{n-1}n-1}\lambda^{r_{n}+r_{n-1}+m_{3}%
}\times\medskip\\
\times\sin\lambda\rho_{n}(\xi_{n}-b)\sin\lambda\rho_{n-1}\left(  \xi_{n}%
-\xi_{n-1}\right)  \sin\lambda\rho_{n-2}\left(  x-\xi_{n-1}\right)  \medskip\\
+o(\lambda^{r_{n}+r_{n-1}+m_{3}}\exp\left\vert \operatorname{Im}%
\lambda\right\vert ((\xi_{n}-b)\rho_{n}+\left(  \xi_{n}-\xi_{n-1}\right)
\rho_{n-1}+\left(  x-\xi_{n-1}\right)  \rho_{n-2})\medskip\\
,\text{ }\deg b_{1}\left(  \lambda\right)  >\deg b_{2}\left(  \lambda\right)
\medskip
\end{array}
\medskip\\%
\begin{array}
[c]{l}%
\gamma_{r_{n}n}\gamma_{r_{n-1}n-1}\lambda^{r_{n}+r_{n-1}+m_{b}}\sin\lambda
\rho_{n-2}\left(  x-\xi_{n-1}\right)  \sin\lambda\rho_{n-1}\left(  \xi_{n}%
-\xi_{n-1}\right)  \times\medskip\\
\times\left(  b_{m_{3}1}\sin\lambda\rho_{n}(\xi_{n}-b)-b_{m_{4}2}\cos
\lambda\rho_{n}(\xi_{n}-b)\right)  \medskip\\
+o(\lambda^{r_{n}+r_{n-1}+m_{b}}\exp\left\vert \operatorname{Im}%
\lambda\right\vert ((\xi_{n}-b)\rho_{n}+\left(  \xi_{n}-\xi_{n-1}\right)
\rho_{n-1}+\left(  x-\xi_{n-1}\right)  \rho_{n-2})\\
,\text{ }\deg b_{1}\left(  \lambda\right)  =\deg b_{2}\left(  \lambda\right)
\medskip
\end{array}
\end{array}
\right.  \smallskip$

$\psi_{n-2,2}(x,\lambda)=\left\{
\begin{array}
[c]{l}%
\begin{array}
[c]{l}%
b_{m_{4}2}\gamma_{r_{n}n}\gamma_{r_{n-1}n-1}\lambda^{r_{n}+r_{n-1}+m_{4}%
}\times\medskip\\
\times\cos\lambda\rho_{n}(\xi_{n}-b)\sin\lambda\rho_{n-1}\left(  \xi_{n}%
-\xi_{n-1}\right)  \cos\lambda\rho_{n-2}\left(  x-\xi_{n-1}\right)  \medskip\\
+o(\lambda^{r_{n}+r_{n-1}+m_{4}}\exp\left\vert \operatorname{Im}%
\lambda\right\vert ((\xi_{n}-b)\rho_{n}+\left(  \xi_{n}-\xi_{n-1}\right)
\rho_{n-1}+\left(  x-\xi_{n-1}\right)  \rho_{n-2})\medskip\\
,\text{ }\deg b_{2}\left(  \lambda\right)  >\deg b_{1}\left(  \lambda\right)
\medskip
\end{array}
\\%
\begin{array}
[c]{l}%
-b_{m_{3}1}\gamma_{r_{n}n}\gamma_{r_{n-1}n-1}\lambda^{r_{n}+r_{n-1}+m_{3}%
}\times\medskip\\
\times\sin\lambda\rho_{n}(\xi_{n}-b)\sin\lambda\rho_{n-1}\left(  \xi_{n}%
-\xi_{n-1}\right)  \sin\lambda\rho_{n-2}\left(  x-\xi_{n-1}\right)  \medskip\\
+o(\lambda^{r_{n}+r_{n-1}+m_{3}}\exp\left\vert \operatorname{Im}%
\lambda\right\vert ((\xi_{n}-b)\rho_{n}+\left(  \xi_{n}-\xi_{n-1}\right)
\rho_{n-1}+\left(  x-\xi_{n-1}\right)  \rho_{n-2})\medskip\\
,\text{ }\deg b_{1}\left(  \lambda\right)  >\deg b_{2}\left(  \lambda\right)
\medskip
\end{array}
\\%
\begin{array}
[c]{l}%
\gamma_{r_{n}n}\gamma_{r_{n-1}n-1}\lambda^{r_{n}+r_{n-1}+m_{b}}\cos\lambda
\rho_{n-2}\left(  x-\xi_{n-1}\right)  \sin\lambda\rho_{n-1}\left(  \xi_{n}%
-\xi_{n-1}\right)  \times\medskip\\
\times\left(  -b_{m_{3}1}\sin\lambda\rho_{n}(\xi_{n}-b)+b_{m_{4}2}\cos
\lambda\rho_{n}(\xi_{n}-b)\right)  \medskip\\
+o(\lambda^{r_{n}+r_{n-1}+m_{b}}\exp\left\vert \operatorname{Im}%
\lambda\right\vert ((\xi_{n}-b)\rho_{n}+\left(  \xi_{n}-\xi_{n-1}\right)
\rho_{n-1}+\left(  x-\xi_{n-1}\right)  \rho_{2})\medskip\\
,\text{ }\deg b_{1}\left(  \lambda\right)  =\deg b_{2}\left(  \lambda\right)
\medskip
\end{array}
\end{array}
\right.  \smallskip$

$,\ldots,$\newline$\psi_{11}(x,\lambda)=\left\{
\begin{array}
[c]{l}%
\begin{array}
[c]{l}%
\left(  -1\right)  ^{n}b_{m_{4}2}\lambda^{r_{2}+\ldots+r_{n}+m_{4}}%
\times\left(
{\displaystyle\prod\limits_{i=2}^{n}}
\gamma_{r_{i}i}\right)  \times\medskip\\
\times\cos\lambda\rho_{n}(\xi_{n}-b)\sin\lambda\rho_{1}\left(  x-\xi
_{2}\right)  \times\left(
{\displaystyle\prod\limits_{i=3}^{n}}
\sin\lambda\rho_{i-1}\left(  \xi_{i}-\xi_{i-1}\right)  \right)  \medskip\\
+o(\lambda^{r_{2}+\ldots+r_{n}+m_{4}}\exp\left\vert \operatorname{Im}%
\lambda\right\vert (\left(  \xi_{2}-\xi_{1}\right)  \rho_{1}+%
{\displaystyle\sum\limits_{i=3}^{n}}
\left(  \xi_{i}-\xi_{i-1}\right)  \rho_{i-1}),\text{ }\deg b_{2}\left(
\lambda\right)  >\deg b_{1}\left(  \lambda\right)  \medskip
\end{array}
\\%
\begin{array}
[c]{l}%
\left(  -1\right)  ^{n-1}b_{m_{3}1}\lambda^{r_{2}+\ldots+r_{n}+m_{3}}%
\times\left(
{\displaystyle\prod\limits_{i=2}^{n}}
\gamma_{r_{i}i}\right)  \times\medskip\\
\times\sin\lambda\rho_{n}(\xi_{n}-b)\sin\lambda\rho_{1}\left(  x-\xi
_{2}\right)  \times\left(
{\displaystyle\prod\limits_{i=3}^{n}}
\sin\lambda\rho_{i-1}\left(  \xi_{i}-\xi_{i-1}\right)  \right)  \medskip\\
+o(\lambda^{r_{2}+\ldots+r_{n}+m_{3}}\exp\left\vert \operatorname{Im}%
\lambda\right\vert (\left(  \xi_{2}-\xi_{1}\right)  \rho_{1}+%
{\displaystyle\sum\limits_{i=3}^{n}}
\left(  \xi_{i}-\xi_{i-1}\right)  \rho_{i-1}),\text{ }\deg b_{1}\left(
\lambda\right)  >\deg b_{2}\left(  \lambda\right)  \medskip
\end{array}
\\%
\begin{array}
[c]{l}%
\lambda^{r_{2}+\ldots+r_{n}+m_{b}}\times\left(
{\displaystyle\prod\limits_{i=2}^{n}}
\gamma_{r_{i}i}\right)  \times\sin\lambda\rho_{1}\left(  x-\xi_{2}\right)
\times\medskip\\
\times\left(
{\displaystyle\prod\limits_{i=3}^{n}}
\sin\lambda\rho_{i-1}\left(  \xi_{i}-\xi_{i-1}\right)  \right)  \left\{
\begin{array}
[c]{c}%
\left(  -1\right)  ^{n-1}b_{m_{3}1}\sin\lambda\rho_{n}(\xi_{n}-b)\medskip\\
+\left(  -1\right)  ^{n}b_{m_{4}2}\cos\lambda\rho_{n}(\xi_{n}-b)
\end{array}
\right\}  \medskip\\
+o(\lambda^{r_{2}+\ldots+r_{n}+m_{b}}\exp\left\vert \operatorname{Im}%
\lambda\right\vert (\left(  \xi_{2}-\xi_{1}\right)  \rho_{1}+%
{\displaystyle\sum\limits_{i=3}^{n}}
\left(  \xi_{i}-\xi_{i-1}\right)  \rho_{i-1}),\text{ }\deg b_{1}\left(
\lambda\right)  =\deg b_{2}\left(  \lambda\right)  \medskip
\end{array}
\smallskip
\end{array}
\right.  \smallskip$

$\psi_{12}(x,\lambda)=\left\{
\begin{array}
[c]{l}%
\begin{array}
[c]{l}%
\left(  -1\right)  ^{n-1}b_{m_{4}2}\lambda^{r_{2}+\ldots+r_{n}+m_{4}}%
\times\left(
{\displaystyle\prod\limits_{i=2}^{n}}
\gamma_{r_{i}i}\right)  \times\medskip\\
\times\cos\lambda\rho_{n}(\xi_{n}-b)\cos\lambda\rho_{1}\left(  x-\xi
_{2}\right)  \times\left(
{\displaystyle\prod\limits_{i=3}^{n}}
\sin\lambda\rho_{i-1}\left(  \xi_{i}-\xi_{i-1}\right)  \right)  \medskip\\
+o(\lambda^{r_{2}+\ldots+r_{n}+m_{4}}\exp\left\vert \operatorname{Im}%
\lambda\right\vert (\left(  \xi_{2}-\xi_{1}\right)  \rho_{1}+%
{\displaystyle\sum\limits_{i=3}^{n}}
\left(  \xi_{i}-\xi_{i-1}\right)  \rho_{i-1}),\text{ }\deg b_{2}\left(
\lambda\right)  >\deg b_{1}\left(  \lambda\right)  \medskip
\end{array}
\\%
\begin{array}
[c]{l}%
\left(  -1\right)  ^{n}b_{m_{3}1}\lambda^{r_{2}+\ldots+r_{n}+m_{3}}%
\times\left(
{\displaystyle\prod\limits_{i=2}^{n}}
\gamma_{r_{i}i}\right)  \times\medskip\\
\times\sin\lambda\rho_{n}(\xi_{n}-b)\cos\lambda\rho_{1}\left(  x-\xi
_{2}\right)  \times\left(
{\displaystyle\prod\limits_{i=3}^{n}}
\sin\lambda\rho_{i-1}\left(  \xi_{i}-\xi_{i-1}\right)  \right)  \medskip\\
+o(\lambda^{r_{2}+\ldots+r_{n}+m_{3}}\exp\left\vert \operatorname{Im}%
\lambda\right\vert (\left(  \xi_{2}-\xi_{1}\right)  \rho_{1}+%
{\displaystyle\sum\limits_{i=3}^{n}}
\left(  \xi_{i}-\xi_{i-1}\right)  \rho_{i-1}),\text{ }\deg b_{1}\left(
\lambda\right)  >\deg b_{2}\left(  \lambda\right)  \medskip
\end{array}
\\%
\begin{array}
[c]{l}%
\lambda^{r_{2}+\ldots+r_{n}+m_{b}}\times\left(
{\displaystyle\prod\limits_{i=2}^{n}}
\gamma_{r_{i}i}\right)  \times\cos\lambda\rho_{1}\left(  x-\xi_{2}\right)
\times\medskip\\
\times\left(
{\displaystyle\prod\limits_{i=3}^{n}}
\sin\lambda\rho_{i-1}\left(  \xi_{i}-\xi_{i-1}\right)  \right)  \left\{
\begin{array}
[c]{c}%
\left(  -1\right)  ^{n}b_{m_{3}1}\sin\lambda\rho_{n}(\xi_{n}-b)\medskip\\
+\left(  -1\right)  ^{n-1}b_{m_{4}2}\cos\lambda\rho_{n}(\xi_{n}-b)\medskip
\end{array}
\right\}  \medskip\\
+o(\lambda^{r_{2}+\ldots+r_{n}+m_{b}}\exp\left\vert \operatorname{Im}%
\lambda\right\vert (\left(  \xi_{2}-\xi_{1}\right)  \rho_{1}+%
{\displaystyle\sum\limits_{i=3}^{n}}
\left(  \xi_{i}-\xi_{i-1}\right)  \rho_{i-1}),\text{ }\deg b_{1}\left(
\lambda\right)  =\deg b_{2}\left(  \lambda\right)  \medskip
\end{array}
\end{array}
\right.  $\smallskip

$\bigskip$

$\psi_{01}(x,\lambda)=\left\{
\begin{array}
[c]{l}%
\begin{array}
[c]{l}%
\left(  -1\right)  ^{n-1}b_{m_{4}2}\lambda^{A+m_{4}}\times\left(
{\displaystyle\prod\limits_{i=1}^{n}}
\gamma_{r_{i}i}\right)  \times\medskip\\
\times\cos\lambda\rho_{n}(\xi_{n}-b)\sin\lambda\rho_{0}\left(  x-\xi
_{1}\right)  \times\left(
{\displaystyle\prod\limits_{i=2}^{n}}
\sin\lambda\rho_{i-1}\left(  \xi_{i}-\xi_{i-1}\right)  \right)  \medskip\\
+o(\lambda^{A+m_{4}}\exp\left\vert \operatorname{Im}\lambda\right\vert
(\left(  \xi_{1}-a\right)  \rho_{0}+%
{\displaystyle\sum\limits_{i=2}^{n}}
\left(  \xi_{i}-\xi_{i-1}\right)  \rho_{i-1}),\text{ }\deg b_{2}\left(
\lambda\right)  >\deg b_{1}\left(  \lambda\right)  \medskip
\end{array}
\\%
\begin{array}
[c]{l}%
\left(  -1\right)  ^{n}b_{m_{3}1}\lambda^{A+m_{3}}\times\left(
{\displaystyle\prod\limits_{i=1}^{n}}
\gamma_{r_{i}i}\right)  \times\medskip\\
\times\sin\lambda\rho_{n}(\xi_{n}-b)\sin\lambda\rho_{0}\left(  x-\xi
_{1}\right)  \times\left(
{\displaystyle\prod\limits_{i=2}^{n}}
\sin\lambda\rho_{i-1}\left(  \xi_{i}-\xi_{i-1}\right)  \right)  \medskip\\
+o(\lambda^{A+m_{3}}\exp\left\vert \operatorname{Im}\lambda\right\vert
(\left(  \xi_{1}-a\right)  \rho_{0}+%
{\displaystyle\sum\limits_{i=2}^{n}}
\left(  \xi_{i}-\xi_{i-1}\right)  \rho_{i-1}),\text{ }\deg b_{1}\left(
\lambda\right)  >\deg b_{2}\left(  \lambda\right)  \medskip
\end{array}
\\%
\begin{array}
[c]{l}%
\lambda^{A+m_{b}}\times\left(
{\displaystyle\prod\limits_{i=1}^{n}}
\gamma_{r_{i}i}\right)  \times\sin\lambda\rho_{0}\left(  x-\xi_{1}\right)
\times\medskip\\
\times\left(
{\displaystyle\prod\limits_{i=2}^{n}}
\sin\lambda\rho_{i-1}\left(  \xi_{i}-\xi_{i-1}\right)  \right)  \left\{
\begin{array}
[c]{c}%
\left(  -1\right)  ^{n}b_{m_{3}1}\sin\lambda\rho_{n}(\xi_{n}-b)\medskip\\
+\left(  -1\right)  ^{n-1}b_{m_{4}2}\cos\lambda\rho_{n}(\xi_{n}-b)\medskip
\end{array}
\right\}  \medskip\\
+o(\lambda^{A+m_{b}}\exp\left\vert \operatorname{Im}\lambda\right\vert
(\left(  \xi_{1}-a\right)  \rho_{0}+%
{\displaystyle\sum\limits_{i=2}^{n}}
\left(  \xi_{i}-\xi_{i-1}\right)  \rho_{i-1}),\text{ }\deg b_{1}\left(
\lambda\right)  =\deg b_{2}\left(  \lambda\right)  \medskip
\end{array}
\end{array}
\right.  \smallskip$

$\psi_{02}(x,\lambda)=\left\{
\begin{array}
[c]{l}%
\begin{array}
[c]{l}%
\left(  -1\right)  ^{n}b_{m_{4}2}\lambda^{A+m_{4}}\times\left(
{\displaystyle\prod\limits_{i=1}^{n}}
\gamma_{r_{i}i}\right)  \times\medskip\\
\times\cos\lambda\rho_{n}(\xi_{n}-b)\cos\lambda\rho_{0}\left(  x-\xi
_{1}\right)  \times\left(
{\displaystyle\prod\limits_{i=2}^{n}}
\sin\lambda\rho_{i-1}\left(  \xi_{i}-\xi_{i-1}\right)  \right)  \medskip\\
+o(\lambda^{A+m_{4}}\exp\left\vert \operatorname{Im}\lambda\right\vert
(\left(  \xi_{1}-a\right)  \rho_{0}+%
{\displaystyle\sum\limits_{i=2}^{n}}
\left(  \xi_{i}-\xi_{i-1}\right)  \rho_{i-1}),\text{ }\deg b_{2}\left(
\lambda\right)  >\deg b_{1}\left(  \lambda\right)  \medskip
\end{array}
\\%
\begin{array}
[c]{l}%
\left(  -1\right)  ^{n-1}b_{m_{3}1}\lambda^{A+m_{3}}\times\left(
{\displaystyle\prod\limits_{i=1}^{n}}
\gamma_{r_{i}i}\right)  \times\medskip\\
\times\sin\lambda\rho_{n}(\xi_{n}-b)\cos\lambda\rho_{0}\left(  x-\xi
_{1}\right)  \times\left(
{\displaystyle\prod\limits_{i=2}^{n}}
\sin\lambda\rho_{i-1}\left(  \xi_{i}-\xi_{i-1}\right)  \right)  \medskip\\
+o(\lambda^{A+m_{3}}\exp\left\vert \operatorname{Im}\lambda\right\vert
(\left(  \xi_{1}-a\right)  \rho_{0}+%
{\displaystyle\sum\limits_{i=2}^{n}}
\left(  \xi_{i}-\xi_{i-1}\right)  \rho_{i-1}),\text{ }\deg b_{1}\left(
\lambda\right)  >\deg b_{2}\left(  \lambda\right)  \medskip
\end{array}
\\%
\begin{array}
[c]{l}%
\lambda^{A+m_{b}}\times\left(
{\displaystyle\prod\limits_{i=1}^{n}}
\gamma_{r_{i}i}\right)  \times\cos\lambda\rho_{0}\left(  x-\xi_{1}\right)
\times\medskip\\
\times\left(
{\displaystyle\prod\limits_{i=2}^{n}}
\sin\lambda\rho_{i-1}\left(  \xi_{i}-\xi_{i-1}\right)  \right)  \left\{
\begin{array}
[c]{c}%
\left(  -1\right)  ^{n-1}b_{m_{3}1}\sin\lambda\rho_{n}(\xi_{n}-b)\medskip\\
+\left(  -1\right)  ^{n}b_{m_{4}2}\cos\lambda\rho_{n}(\xi_{n}-b)\medskip
\end{array}
\right\}  \medskip\\
+o(\lambda^{A+m_{b}}\exp\left\vert \operatorname{Im}\lambda\right\vert
(\left(  \xi_{1}-a\right)  \rho_{0}+%
{\displaystyle\sum\limits_{i=2}^{n}}
\left(  \xi_{i}-\xi_{i-1}\right)  \rho_{i-1}),\text{ }\deg b_{1}\left(
\lambda\right)  =\deg b_{2}\left(  \lambda\right)  \text{.}\medskip
\end{array}
\end{array}
\right.  $

\bigskip

The function%
\begin{equation}%
\begin{array}
[c]{ll}%
\Delta(\lambda)=W\left(  \psi,\varphi\right)  & =\psi\left(  x,\lambda\right)
\varphi^{\prime}\left(  x,\lambda\right)  -\psi^{\prime}\left(  x,\lambda
\right)  \varphi\left(  x,\lambda\right) \\
& =\psi_{1}(b,\lambda)\varphi_{2}(b,\lambda)-\psi_{2}(b,\lambda)\varphi
_{1}(b,\lambda)\\
& =\psi_{1}(a,\lambda)\varphi_{2}(a,\lambda)-\psi_{2}(a,\lambda)\varphi
_{1}(a,\lambda)\\
& =b_{2}\left(  \lambda\right)  \varphi_{2}(b,\lambda)-b_{1}\left(
\lambda\right)  \varphi_{1}(b,\lambda)\\
& =a_{1}\left(  \lambda\right)  \psi_{1}(a,\lambda)-a_{2}\left(
\lambda\right)  \psi_{2}(a,\lambda)
\end{array}
\tag{9}%
\end{equation}
is called characteristic function of the problem $L$. Besides, $\Delta
(\lambda)$ is an entire function in $\lambda$ and it is obvious that zeros of
$\Delta(\lambda)$ coincide with the eigenvalues of problem $L$.

\smallskip\textbf{Lemma 3: }The following asymptotic behaviours hold
$\left\vert \lambda\right\vert \rightarrow\infty$:$\smallskip$\newline for
$\deg a_{2}\left(  \lambda\right)  >\deg a_{1}\left(  \lambda\right)
;\medskip$

$\Delta(\lambda)=\left\{
\begin{array}
[c]{l}%
\begin{array}
[c]{l}%
\left(  -1\right)  ^{n+1}a_{m_{2}2}b_{m_{4}2}\lambda^{A+m_{2}+m_{4}}%
\times\medskip\\
\times\left(
{\displaystyle\prod\limits_{i=1}^{n}}
\gamma_{r_{i}i}\right)  \cos\lambda\rho_{0}(\xi_{1}-a)\cos\lambda\rho
_{n}(x-\xi_{n})\times\left(
{\displaystyle\prod\limits_{i=2}^{n}}
\sin\lambda\rho_{i-1}\left(  \xi_{i}-\xi_{i-1}\right)  \right)  \medskip\\
+o(\lambda^{A+m_{2}+m_{4}}\exp\left\vert \operatorname{Im}\lambda\right\vert
(\left(  b-\xi_{n}\right)  \rho_{n}+%
{\displaystyle\sum\limits_{i=1}^{n}}
\left(  \xi_{i}-\xi_{i-1}\right)  \rho_{i-1}),\text{ }\deg b_{2}\left(
\lambda\right)  >\deg b_{1}\left(  \lambda\right)  \medskip
\end{array}
\smallskip\\%
\begin{array}
[c]{l}%
\left(  -1\right)  ^{n+1}a_{m_{2}2}b_{m_{3}1}\lambda^{A+m_{2}+m_{3}}%
\times\medskip\\
\times\left(
{\displaystyle\prod\limits_{i=1}^{n}}
\gamma_{r_{i}i}\right)  \cos\lambda\rho_{0}(\xi_{1}-a)\sin\lambda\rho
_{n}(x-\xi_{n})\times\left(
{\displaystyle\prod\limits_{i=2}^{n}}
\sin\lambda\rho_{i-1}\left(  \xi_{i}-\xi_{i-1}\right)  \right)  \medskip\\
+o(\lambda^{A+m_{2}+m_{3}}\exp\left\vert \operatorname{Im}\lambda\right\vert
(\left(  b-\xi_{n}\right)  \rho_{n}+%
{\displaystyle\sum\limits_{i=1}^{n}}
\left(  \xi_{i}-\xi_{i-1}\right)  \rho_{i-1}),\text{ }\deg b_{1}\left(
\lambda\right)  >\deg b_{2}\left(  \lambda\right)  \medskip
\end{array}
\smallskip\\%
\begin{array}
[c]{l}%
\left(  -1\right)  ^{n+1}\lambda^{A+m_{2}+m_{b}}a_{m_{2}2}\times\left(
{\displaystyle\prod\limits_{i=1}^{n}}
\gamma_{r_{i}i}\right)  \times\medskip\\
\times\cos\lambda\rho_{0}(\xi_{1}-a)\times\left(
{\displaystyle\prod\limits_{i=2}^{n}}
\sin\lambda\rho_{i-1}\left(  \xi_{i}-\xi_{i-1}\right)  \right)  \left\{
\begin{array}
[c]{c}%
b_{m_{4}2}\cos\lambda\rho_{n}(x-\xi_{n})\medskip\\
+b_{m_{3}1}\sin\lambda\rho_{n}(x-\xi_{n})\medskip
\end{array}
\right\}  \medskip\\
+o(\lambda^{A+m_{2}+m_{b}}\exp\left\vert \operatorname{Im}\lambda\right\vert
(\left(  b-\xi_{n}\right)  \rho_{n}+%
{\displaystyle\sum\limits_{i=1}^{n}}
\left(  \xi_{i}-\xi_{i-1}\right)  \rho_{i-1}),\text{ }\deg b_{1}\left(
\lambda\right)  =\deg b_{2}\left(  \lambda\right)  \medskip
\end{array}
\smallskip
\end{array}
\right.  \smallskip$

for $\deg a_{1}\left(  \lambda\right)  >\deg a_{2}\left(  \lambda\right)
;\medskip$

$\Delta(\lambda)=\left\{
\begin{array}
[c]{l}%
\begin{array}
[c]{l}%
\left(  -1\right)  ^{n}a_{m_{1}1}b_{m_{4}2}\lambda^{A+m_{1}+m_{4}}%
\times\medskip\\
\times\left(
{\displaystyle\prod\limits_{i=1}^{n}}
\gamma_{r_{i}i}\right)  \sin\lambda\rho_{0}(\xi_{1}-a)\cos\lambda\rho
_{n}(x-\xi_{n})\times\left(
{\displaystyle\prod\limits_{i=2}^{n}}
\sin\lambda\rho_{i-1}\left(  \xi_{i}-\xi_{i-1}\right)  \right)  \medskip\\
+o(\lambda^{A+m_{1}+m_{4}}\exp\left\vert \operatorname{Im}\lambda\right\vert
(\left(  b-\xi_{n}\right)  \rho_{n}+%
{\displaystyle\sum\limits_{i=1}^{n}}
\left(  \xi_{i}-\xi_{i-1}\right)  \rho_{i-1}),\text{ }\deg b_{2}\left(
\lambda\right)  >\deg b_{1}\left(  \lambda\right)  \medskip
\end{array}
\smallskip\\%
\begin{array}
[c]{l}%
\left(  -1\right)  ^{n}a_{m_{1}1}b_{m_{3}1}\lambda^{A+m_{1}+m_{3}}%
\times\medskip\\
\times\left(
{\displaystyle\prod\limits_{i=1}^{n}}
\gamma_{r_{i}i}\right)  \sin\lambda\rho_{0}(\xi_{1}-a)\sin\lambda\rho
_{n}(x-\xi_{n})\times\left(
{\displaystyle\prod\limits_{i=2}^{n}}
\sin\lambda\rho_{i-1}\left(  \xi_{i}-\xi_{i-1}\right)  \right)  \medskip\\
+o(\lambda^{A+m_{1}+m_{3}}\exp\left\vert \operatorname{Im}\lambda\right\vert
(\left(  b-\xi_{n}\right)  \rho_{n}+%
{\displaystyle\sum\limits_{i=1}^{n}}
\left(  \xi_{i}-\xi_{i-1}\right)  \rho_{i-1}),\text{ }\deg b_{1}\left(
\lambda\right)  >\deg b_{2}\left(  \lambda\right)  \medskip
\end{array}
\smallskip\\%
\begin{array}
[c]{l}%
\left(  -1\right)  ^{n}\lambda^{A+m_{1}+m_{b}}a_{m_{1}1}\times\left(
{\displaystyle\prod\limits_{i=1}^{n}}
\gamma_{r_{i}i}\right)  \times\medskip\\
\times\sin\lambda\rho_{0}(\xi_{1}-a)\times\left(
{\displaystyle\prod\limits_{i=2}^{n}}
\sin\lambda\rho_{i-1}\left(  \xi_{i}-\xi_{i-1}\right)  \right)  \left\{
\begin{array}
[c]{c}%
b_{m_{4}2}\cos\lambda\rho_{n}(x-\xi_{n})\\
+b_{m_{3}1}\sin\lambda\rho_{n}(x-\xi_{n})
\end{array}
\right\}  \medskip\\
+o(\lambda^{A+m_{1}+m_{b}}\exp\left\vert \operatorname{Im}\lambda\right\vert
(\left(  b-\xi_{n}\right)  \rho_{n}+%
{\displaystyle\sum\limits_{i=1}^{n}}
\left(  \xi_{i}-\xi_{i-1}\right)  \rho_{i-1}),\text{ }\deg b_{1}\left(
\lambda\right)  =\deg b_{2}\left(  \lambda\right)  \medskip
\end{array}
\smallskip
\end{array}
\right.  \smallskip$

\bigskip

\bigskip for $\deg a_{1}\left(  \lambda\right)  =\deg a_{2}\left(
\lambda\right)  ;\medskip$

$\Delta(\lambda)=\left\{
\begin{array}
[c]{l}%
\begin{array}
[c]{l}%
\lambda^{A+m+m_{4}}\times\left(
{\displaystyle\prod\limits_{i=1}^{n}}
\gamma_{r_{i}i}\right)  \times\left(
{\displaystyle\prod\limits_{i=2}^{n}}
\sin\lambda\rho_{i-1}\left(  \xi_{i}-\xi_{i-1}\right)  \right)  \times
\medskip\\
\times\left(
\begin{array}
[c]{c}%
\left(  -1\right)  ^{n+1}a_{m_{2}2}b_{m_{4}2}\cos\lambda\rho_{0}(\xi
_{1}-a)\cos\lambda\rho_{n}(x-\xi_{n})\medskip\\
+\left(  -1\right)  ^{n}a_{m_{1}1}b_{m_{4}2}\sin\lambda\rho_{0}(\xi_{1}%
-a)\cos\lambda\rho_{n}(x-\xi_{n})\medskip
\end{array}
\right)  \medskip\\
+o(\lambda^{A+m+m_{4}}\exp\left\vert \operatorname{Im}\lambda\right\vert
(\left(  b-\xi_{n}\right)  \rho_{n}+%
{\displaystyle\sum\limits_{i=1}^{n}}
\left(  \xi_{i}-\xi_{i-1}\right)  \rho_{i-1}),\text{ }\deg b_{2}\left(
\lambda\right)  >\deg b_{1}\left(  \lambda\right)  \medskip
\end{array}
\smallskip\\%
\begin{array}
[c]{l}%
\lambda^{A+m+m_{3}}\times\left(
{\displaystyle\prod\limits_{i=1}^{n}}
\gamma_{r_{i}i}\right)  \times\left(
{\displaystyle\prod\limits_{i=2}^{n}}
\sin\lambda\rho_{i-1}\left(  \xi_{i}-\xi_{i-1}\right)  \right)  \times
\medskip\\
\times\left(
\begin{array}
[c]{c}%
\left(  -1\right)  ^{n+1}a_{m_{2}2}b_{m_{3}1}\cos\lambda\rho_{0}(\xi
_{1}-a)\sin\lambda\rho_{n}(x-\xi_{n})\medskip\\
+\left(  -1\right)  ^{n+1}a_{m_{1}1}b_{m_{3}1}\sin\lambda\rho_{0}(\xi
_{1}-a)\sin\lambda\rho_{n}(x-\xi_{n})\medskip
\end{array}
\right)  \medskip\\
+o(\lambda^{A+m+m_{3}}\exp\left\vert \operatorname{Im}\lambda\right\vert
(\left(  b-\xi_{n}\right)  \rho_{n}+%
{\displaystyle\sum\limits_{i=1}^{n}}
\left(  \xi_{i}-\xi_{i-1}\right)  \rho_{i-1}),\text{ }\deg b_{1}\left(
\lambda\right)  >\deg b_{2}\left(  \lambda\right)  \medskip
\end{array}
\smallskip\\%
\begin{array}
[c]{l}%
\lambda^{A+m+m_{b}}\times\left(
{\displaystyle\prod\limits_{i=1}^{n}}
\gamma_{r_{i}i}\right)  \times\left(
{\displaystyle\prod\limits_{i=2}^{n}}
\sin\lambda\rho_{i-1}\left(  \xi_{i}-\xi_{i-1}\right)  \right)  \times
\medskip\\
\times\left\{
\begin{array}
[c]{c}%
\left(  -1\right)  ^{n+1}a_{m_{2}2}\cos\lambda\rho_{0}(\xi_{1}-a)\left(
\begin{array}
[c]{c}%
b_{m_{4}2}\cos\lambda\rho_{n}(x-\xi_{n})\medskip\\
+b_{m_{3}1}\sin\lambda\rho_{n}(x-\xi_{n})\medskip
\end{array}
\right)  \medskip\\
+\left(  -1\right)  ^{n}a_{m_{1}1}\sin\lambda\rho_{0}(\xi_{1}-a)\left(
\begin{array}
[c]{c}%
b_{m_{4}2}\cos\lambda\rho_{n}(x-\xi_{n})\medskip\\
+b_{m_{3}1}\sin\lambda\rho_{n}(x-\xi_{n})\medskip
\end{array}
\right)  \medskip
\end{array}
\right\}  \medskip\\
+o(\lambda^{A+m+m_{b}}\exp\left\vert \operatorname{Im}\lambda\right\vert
(\left(  b-\xi_{n}\right)  \rho_{n}+%
{\displaystyle\sum\limits_{i=1}^{n}}
\left(  \xi_{i}-\xi_{i-1}\right)  \rho_{i-1}),\text{ }\deg b_{1}\left(
\lambda\right)  =\deg b_{2}\left(  \lambda\right)  \medskip
\end{array}
\smallskip
\end{array}
\right.  \smallskip$

where $A=r_{1}+r_{2}+\cdots+r_{n}.$

\bigskip

\section{\textbf{Inverse\ Problems}}

In this section,we study the inverse problems for the reconstruction of
boundary value problem (1)-(5) by Weyl function and spectral data.

We consider the boundary value problem $\tilde{L}$ which has the same form
with $L$ but with different coefficients $\tilde{\Omega}(x),$ $\tilde{f}%
_{1}\left(  \lambda\right)  ,$ $\tilde{f}_{2}\left(  \lambda\right)
,\tilde{\varepsilon}_{i},\tilde{\theta}_{i},\tilde{\gamma}_{i}$ such that
\newline$\tilde{\Omega}(x)=\left(
\begin{array}
[c]{ll}%
\tilde{p}\left(  x\right)  & q(x)\\
q(x) & \tilde{r}\left(  x\right)
\end{array}
\right)  .,f_{1}\left(  \lambda\right)  :=\dfrac{a_{1}\left(  \lambda\right)
}{a_{2}\left(  \lambda\right)  },f_{2}\left(  \lambda\right)  :=\dfrac
{b_{1}\left(  \lambda\right)  }{b_{2}\left(  \lambda\right)  }.$

If a certain symbol $\sigma$ denotes an object related to $L,$ then the symbol
$\tilde{\sigma}$ denotes the corresponding object related to $\tilde{L}%
.$\bigskip

\textbf{Theorem 1} If $M(\lambda)=\tilde{M}(\lambda),$ then $L=\tilde{L},$
i.e., $\Omega(x)=\tilde{\Omega}(x),$ a.e. and $f_{1}\left(  \lambda\right)
=\tilde{f}_{1}\left(  \lambda\right)  ,f_{2}\left(  \lambda\right)  =\tilde
{f}_{2}\left(  \lambda\right)  ,\varepsilon_{i}=\tilde{\varepsilon}_{i}%
,\theta_{i}=\tilde{\theta}_{i},\gamma_{i}(\lambda)=\tilde{\gamma}_{i}\left(
\lambda\right)  ,i=\overline{1,n}.$

\textbf{Proof} We introduce a matrix $P(x,\lambda)=\left[  P_{kj}%
(x,\lambda)\right]  _{k,j=1,2}$ by the formula%
\begin{equation}
P(x,\lambda)\left(
\begin{array}
[c]{cc}%
\tilde{\varphi}_{1} & \tilde{\Phi}_{1}\\
\tilde{\varphi}_{2} & \tilde{\Phi}_{2}%
\end{array}
\right)  =\left(
\begin{array}
[c]{cc}%
\varphi_{1} & \Phi_{1}\\
\varphi_{2} & \Phi_{2}%
\end{array}
\right)  \tag{10}%
\end{equation}

or%
\begin{equation}
\left(
\begin{array}
[c]{cc}%
P_{11}(x,\lambda) & P_{12}(x,\lambda)\\
P_{21}(x,\lambda) & P_{22}(x,\lambda)
\end{array}
\right)  =\left(
\begin{array}
[c]{cc}%
\varphi_{1}\tilde{\Phi}_{2}-\Phi_{1}\tilde{\varphi}_{2} & -\varphi_{1}%
\tilde{\Phi}_{1}+\Phi_{1}\tilde{\varphi}_{1}\\
\varphi_{2}\tilde{\Phi}_{2}-\tilde{\varphi}_{2}\Phi_{2} & -\varphi_{2}%
\tilde{\Phi}_{1}+\tilde{\varphi}_{1}\Phi_{2}%
\end{array}
\right)  \tag{11}%
\end{equation}
where $\Phi(x,\lambda)=\dfrac{\psi(x,\lambda)}{\Delta(\lambda)}$ and $W\left(
\tilde{\Phi},\tilde{\varphi}\right)  =1$. Thus, we find%
\begin{equation}%
\begin{array}
[c]{ll}%
P_{11}(x,\lambda) & =\varphi_{1}(x,\lambda)\dfrac{\tilde{\psi}_{2}(x,\lambda
)}{\tilde{\Delta}(x,\lambda)}-\dfrac{\psi_{1}(x,\lambda)}{\Delta(\lambda
)}\tilde{\varphi}_{2}(x,\lambda)\\
P_{12}(x,\lambda) & =-\varphi_{1}(x,\lambda)\dfrac{\tilde{\psi}_{1}%
(x,\lambda)}{\tilde{\Delta}(\lambda)}+\dfrac{\psi_{1}(x,\lambda)}%
{\Delta(\lambda)}\tilde{\varphi}_{1}(x,\lambda)\\
P_{21}(x,\lambda) & =\varphi_{2}(x,\lambda)\dfrac{\tilde{\psi}_{2}(x,\lambda
)}{\tilde{\Delta}(\lambda)}-\dfrac{\psi_{2}(x,\lambda)}{\Delta(\lambda)}%
\tilde{\varphi}_{2}(x,\lambda)\\
P_{22}(x,\lambda) & =-\varphi_{2}(x,\lambda)\dfrac{\tilde{\psi}_{1}%
(x,\lambda)}{\tilde{\Delta}(\lambda)}+\dfrac{\psi_{2}(x,\lambda)}%
{\Delta(\lambda)}\tilde{\varphi}_{1}(x,\lambda)
\end{array}
\tag{12}%
\end{equation}
\newline We get%
\begin{equation}%
\begin{array}
[c]{ll}%
P_{11}(x,\lambda)= & \varphi_{1}(x,\lambda)\tilde{\phi}_{2}\left(
x,\lambda\right)  -\tilde{\varphi}_{2}(x,\lambda)\phi_{1}\left(
x,\lambda\right)  +\left(  \tilde{M}\left(  \lambda\right)  -M\left(
\lambda\right)  \right)  \varphi_{1}(x,\lambda)\tilde{\varphi}_{2}%
(x,\lambda)\\
P_{12}(x,\lambda)= & -\varphi_{1}(x,\lambda)\tilde{\phi}_{1}\left(
x,\lambda\right)  +\tilde{\varphi}_{1}(x,\lambda)\phi_{1}\left(
x,\lambda\right)  -\left(  \tilde{M}\left(  \lambda\right)  -M\left(
\lambda\right)  \right)  \varphi_{1}(x,\lambda)\tilde{\varphi}_{1}%
(x,\lambda)\\
P_{21}(x,\lambda)= & \varphi_{2}(x,\lambda)\tilde{\phi}_{2}\left(
x,\lambda\right)  -\tilde{\varphi}_{2}(x,\lambda)\phi_{2}\left(
x,\lambda\right)  +\left(  \tilde{M}\left(  \lambda\right)  -M\left(
\lambda\right)  \right)  \varphi_{2}(x,\lambda)\tilde{\varphi}_{2}%
(x,\lambda)\\
P_{22}(x,\lambda)= & -\varphi_{2}(x,\lambda)\tilde{\phi}_{1}\left(
x,\lambda\right)  +\tilde{\varphi}_{1}(x,\lambda)\phi_{2}\left(
x,\lambda\right)  -\left(  \tilde{M}\left(  \lambda\right)  -M\left(
\lambda\right)  \right)  \tilde{\varphi}_{1}(x,\lambda)\varphi_{2}(x,\lambda).
\end{array}
\tag{13}%
\end{equation}

Thus, if $M(\lambda)\equiv\tilde{M}(\lambda)$ then the functions
$P_{ij}(x,\lambda)$ $\left(  i,j=1,2\right)  $ are entire in$\smallskip$
$\lambda$ for each fixed $x$. Moreover, since asymptotic behaviours of
$\varphi_{i}(x,\lambda),$ $\tilde{\varphi}_{i}(x,\lambda),$ $\psi
_{i}(x,\lambda),$ $\tilde{\psi}_{i}(x,\lambda)$ and

for $\deg a_{2}\left(  \lambda\right)  >\deg a_{1}\left(  \lambda\right)
;\smallskip$

$\left\vert \Delta(\lambda)\right\vert \geq C_{\delta}\left\{
\begin{array}
[c]{c}%
\left\vert \lambda\right\vert ^{A+m_{2}+m_{4}}\exp\left\vert \operatorname{Im}%
\lambda\right\vert \left(  \left(  b-\xi_{n}\right)  \rho_{n}+%
{\displaystyle\sum\limits_{i=1}^{n}}
\left(  \xi_{i}-\xi_{i-1}\right)  \rho_{i-1}\right)  ,\text{ }\deg
b_{2}\left(  \lambda\right)  >\deg b_{1}\left(  \lambda\right)  \smallskip\\
\left\vert \lambda\right\vert ^{A+m_{2}+m_{3}}\exp\left\vert \operatorname{Im}%
\lambda\right\vert \left(  \left(  b-\xi_{n}\right)  \rho_{n}+%
{\displaystyle\sum\limits_{i=1}^{n}}
\left(  \xi_{i}-\xi_{i-1}\right)  \rho_{i-1}\right)  ,\text{ }\deg
b_{1}\left(  \lambda\right)  >\deg b_{2}\left(  \lambda\right)  \smallskip\\
\left\vert \lambda\right\vert ^{A+m_{2}+m_{b}}\exp\left\vert \operatorname{Im}%
\lambda\right\vert \left(  \left(  b-\xi_{n}\right)  \rho_{n}+%
{\displaystyle\sum\limits_{i=1}^{n}}
\left(  \xi_{i}-\xi_{i-1}\right)  \rho_{i-1}\right)  ,\text{ }\deg
b_{1}\left(  \lambda\right)  =\deg b_{2}\left(  \lambda\right)  \smallskip
\end{array}
\medskip\right.  $

for $\deg a_{1}\left(  \lambda\right)  >\deg a_{2}\left(  \lambda\right)
;\smallskip$

$\left\vert \Delta(\lambda)\right\vert \geq C_{\delta}\left\{
\begin{array}
[c]{c}%
\left\vert \lambda\right\vert ^{A+m_{1}+m_{4}}\exp\left\vert \operatorname{Im}%
\lambda\right\vert \left(  \left(  b-\xi_{n}\right)  \rho_{n}+%
{\displaystyle\sum\limits_{i=1}^{n}}
\left(  \xi_{i}-\xi_{i-1}\right)  \rho_{i-1}\right)  ,\text{ }\deg
b_{2}\left(  \lambda\right)  >\deg b_{1}\left(  \lambda\right)  \smallskip\\
\left\vert \lambda\right\vert ^{A+m_{1}+m_{3}}\exp\left\vert \operatorname{Im}%
\lambda\right\vert \left(  \left(  b-\xi_{n}\right)  \rho_{n}+%
{\displaystyle\sum\limits_{i=1}^{n}}
\left(  \xi_{i}-\xi_{i-1}\right)  \rho_{i-1}\right)  ,\text{ }\deg
b_{1}\left(  \lambda\right)  >\deg b_{2}\left(  \lambda\right)  \smallskip\\
\left\vert \lambda\right\vert ^{A+m_{1}+m_{b}}\exp\left\vert \operatorname{Im}%
\lambda\right\vert \left(  \left(  b-\xi_{n}\right)  \rho_{n}+%
{\displaystyle\sum\limits_{i=1}^{n}}
\left(  \xi_{i}-\xi_{i-1}\right)  \rho_{i-1}\right)  ,\text{ }\deg
b_{1}\left(  \lambda\right)  =\deg b_{2}\left(  \lambda\right)  \smallskip
\end{array}
\medskip\right.  $

for $\deg a_{1}\left(  \lambda\right)  =\deg a_{2}\left(  \lambda\right)
;\smallskip$

$\left\vert \Delta(\lambda)\right\vert \geq C_{\delta}\left\{
\begin{array}
[c]{c}%
\left\vert \lambda\right\vert ^{A+m+m_{4}}\exp\left\vert \operatorname{Im}%
\lambda\right\vert \left(  \left(  b-\xi_{n}\right)  \rho_{n}+%
{\displaystyle\sum\limits_{i=1}^{n}}
\left(  \xi_{i}-\xi_{i-1}\right)  \rho_{i-1}\right)  ,\text{ }\deg
b_{2}\left(  \lambda\right)  >\deg b_{1}\left(  \lambda\right)  \smallskip\\
\left\vert \lambda\right\vert ^{A+m+m_{3}}\exp\left\vert \operatorname{Im}%
\lambda\right\vert \left(  \left(  b-\xi_{n}\right)  \rho_{n}+%
{\displaystyle\sum\limits_{i=1}^{n}}
\left(  \xi_{i}-\xi_{i-1}\right)  \rho_{i-1}\right)  ,\text{ }\deg
b_{1}\left(  \lambda\right)  >\deg b_{2}\left(  \lambda\right)  \smallskip\\
\left\vert \lambda\right\vert ^{A+mm_{b}}\exp\left\vert \operatorname{Im}%
\lambda\right\vert \left(  \left(  b-\xi_{n}\right)  \rho_{n}+%
{\displaystyle\sum\limits_{i=1}^{n}}
\left(  \xi_{i}-\xi_{i-1}\right)  \rho_{i-1}\right)  ,\text{ }\deg
b_{1}\left(  \lambda\right)  =\deg b_{2}\left(  \lambda\right)  \smallskip
\end{array}
\medskip\right.  $

in $F_{\delta}\cap\tilde{F}_{\delta}$, we$\smallskip$ can easily see that
functions $P_{ij}(x,\lambda)$ are bounded with respect to $\lambda$. As a
result$\smallskip$, these functions don't depend on $\lambda$.

Here, we denote $\tilde{F}_{\delta}=\left\{  \lambda:\left\vert \lambda
-\tilde{\lambda}_{n}\right\vert \geq\delta,\text{ }n=0,\pm1,\pm2,\ldots
\right\}  $ where $n$ is sufficiently small number, $\tilde{\lambda}_{n}$ are
eigenvalues of the problem $\tilde{L}$.

From (12),$\smallskip$

$P_{11}(x,\lambda)-1=\dfrac{\tilde{\psi}_{2}(x,\lambda)\left(  \varphi
_{1}(x,\lambda)-\tilde{\varphi}_{1}(x,\lambda)\right)  }{\tilde{\Delta}\left(
\lambda\right)  }-\tilde{\varphi}_{2}(x,\lambda)\left(  \dfrac{\psi
_{1}(x,\lambda)}{\Delta(\lambda)}-\dfrac{\tilde{\psi}_{1}(x,\lambda)}%
{\tilde{\Delta}(\lambda)}\right)  $

$P_{12}(x,\lambda)=\dfrac{\psi_{1}(x,\lambda)\left(  \tilde{\varphi}%
_{1}(x,\lambda)-\varphi_{1}(x,\lambda)\right)  }{\Delta\left(  \lambda\right)
}+\varphi_{1}(x,\lambda)\left(  \dfrac{\psi_{1}(x,\lambda)}{\Delta(\lambda
)}-\dfrac{\tilde{\psi}_{1}(x,\lambda)}{\tilde{\Delta}(\lambda)}\right)  $

$P_{21}(x,\lambda)=\varphi_{2}(x,\lambda)\left(  \dfrac{\tilde{\psi}%
_{2}(x,\lambda)}{\tilde{\Delta}(\lambda)}-\dfrac{\psi_{2}(x,\lambda)}%
{\Delta\left(  \lambda\right)  }\right)  +\psi_{2}(x,\lambda)\left(
\dfrac{\varphi_{2}(x,\lambda)-\tilde{\varphi}_{2}(x,\lambda)}{\Delta(\lambda
)}\right)  $

$P_{22}(x,\lambda)-1=\dfrac{\psi_{2}(x,\lambda)\left(  \tilde{\varphi}%
_{1}(x,\lambda)-\varphi_{1}(x,\lambda)\right)  }{\Delta\left(  \lambda\right)
}+\varphi_{2}(x,\lambda)\left(  \dfrac{\psi_{1}(x,\lambda)}{\Delta(\lambda
)}-\dfrac{\tilde{\psi}_{1}(x,\lambda)}{\tilde{\Delta}(\lambda)}\right)
\smallskip$

It follows from the representations of solutions $\varphi(x,\lambda)$ and
$\psi(x,\lambda)$,$\smallskip$\newline$\underset{\lambda\rightarrow
\infty}{\lim}\dfrac{\tilde{\psi}_{2}(x,\lambda)\left(  \varphi_{1}%
(x,\lambda)-\tilde{\varphi}_{1}(x,\lambda)\right)  }{\tilde{\Delta}\left(
\lambda\right)  }=0$ and\newline$\underset{\lambda\rightarrow\infty}{\lim
}\tilde{\varphi}_{2}(x,\lambda)\left(  \dfrac{\psi_{1}(x,\lambda)}%
{\Delta(\lambda)}-\dfrac{\tilde{\psi}_{1}(x,\lambda)}{\tilde{\Delta}(\lambda
)}\right)  =0$ for all $x$ in $\Lambda$. Thus,\newline$\underset{\lambda
\rightarrow\infty}{\lim}\left(  P_{11}(x,\lambda)-1\right)  =0$ is valid
uniformly with respect to $x$.$\smallskip$ So we have $P_{11}(x,\lambda
)\equiv1$ and similarly $P_{12}(x,\lambda)\equiv0$, $P_{21}(x,\lambda)\equiv
0$, $P_{22}(x,\lambda)\equiv1.\smallskip$

From (10), we obtain $\varphi_{1}(x,\lambda)\equiv\tilde{\varphi}%
_{1}(x,\lambda)$, $\Phi_{1}\equiv\tilde{\Phi}_{1}$, $\varphi_{2}%
(x,\lambda)\equiv\tilde{\varphi}_{2}(x,\lambda)$ and \newline$\Phi_{2}%
\equiv\tilde{\Phi}_{2}$ for all $x$ and $\lambda$. Moreover, from
$\Phi(x,\lambda)=\dfrac{\psi(x,\lambda)}{\Delta(\lambda)}$, we get\newline%
$\dfrac{\psi_{2}\left(  x,\lambda\right)  }{\psi_{1}\left(  x,\lambda\right)
}=\dfrac{\tilde{\psi}_{2}\left(  x,\lambda\right)  }{\tilde{\psi}_{1}\left(
x,\lambda\right)  }$. Hence, $\Omega(x)=\tilde{\Omega}(x)$, i.e., $p\left(
x\right)  =\tilde{p}\left(  x\right)  $, $r\left(  x\right)  =\tilde{r}\left(
x\right)  $ almost everywhere. On the other hand, since $\left(
\begin{array}
[c]{c}%
\varphi_{01}(a,\lambda)\\
\varphi_{02}(a,\lambda)
\end{array}
\right)  =\left(
\begin{array}
[c]{c}%
a_{2}(\lambda)\\
a_{1}(\lambda)
\end{array}
\right)  $, $\left(
\begin{array}
[c]{c}%
\psi_{n1}(b,\lambda)\\
\psi_{n2}(b,\lambda)
\end{array}
\right)  =\left(
\begin{array}
[c]{c}%
b_{2}(\lambda)\\
b_{1}(\lambda)
\end{array}
\right)  $, we get easily that $f_{1}\left(  \lambda\right)  =\tilde{f}%
_{1}\left(  \lambda\right)  ,f_{2}\left(  \lambda\right)  =\tilde{f}%
_{2}\left(  \lambda\right)  ,\varepsilon_{i}=\tilde{\varepsilon}_{i}%
,\theta_{i}=\tilde{\theta}_{i},$\newline$\gamma_{i}(\lambda)=\tilde{\gamma
}_{i}\left(  \lambda\right)  ,i=\overline{1,n}$. Therefore, $L\equiv\tilde{L}%
$.\smallskip

\textbf{Theorem 2 }If $\lambda_{n}=\tilde{\lambda}_{n}$ and $\mu_{n}%
=\tilde{\mu}_{n}$ for all $n$, then $L\equiv\tilde{L},$ i.e., $\Omega
(x)=\tilde{\Omega}(x),$ a.e.,$f_{1}\left(  \lambda\right)  =\tilde{f}%
_{1}\left(  \lambda\right)  ,f_{2}\left(  \lambda\right)  =\tilde{f}%
_{2}\left(  \lambda\right)  ,\gamma_{i}(\lambda)=\tilde{\gamma}_{i}\left(
\lambda\right)  ,\varepsilon_{i}=\tilde{\varepsilon}_{i},\theta_{i}%
=\tilde{\theta}_{i},i=\overline{1,n}.$ Hence, the problem (1)-(5) is uniquely
determined by spectral data $\left\{  \lambda_{n},\mu_{n}\right\}  .$

\textbf{Proof }If $\lambda_{n}=\tilde{\lambda}_{n}$ and $\mu_{n}=\tilde{\mu
}_{n}$ for all $n$, then $M(\lambda)=\tilde{M}(\lambda)$. Therefore, we get
$L=\tilde{L}$ by Theorem 1$.\medskip$

Let us consider the boundary-value problem $L_{1}$ with the condition\newline%
$y_{1}(a)=0$ instead of the condition (2) in $L$. Let $\left\{  \tau
_{n}\right\}  _{n\in%
\mathbb{Z}
}$ be the eigenvalues of the problem $L_{1}$. It is clear that $\tau_{n}$ are
zeros of\newline$\Delta_{1}(\tau):=\psi_{n1}(a,\tau)$ which is characteristic
function of $L_{1}$.\smallskip

\textbf{Theorem 3 }If $\lambda_{n}=\tilde{\lambda}_{n}$ and $\tau_{n}%
=\tilde{\tau}_{n}$ for all $n$, then\newline$L(\Omega,f_{2}\left(
\lambda\right)  )=L(\tilde{\Omega},\tilde{f}_{2}\left(  \lambda\right)  ).$

Hence, the problem $L$ is uniquely determined by the sequences $\left\{
\lambda_{n}\right\}  $ and $\left\{  \tau_{n}\right\}  $ except coefficients
$f_{1}\left(  \lambda\right)  .$

\textbf{Proof }Since the characteristic functions $\Delta(\lambda)$ and
$\Delta_{1}(\tau)$ are entire of order $1$, functions $\Delta(\lambda)$ and
$\Delta_{1}(\tau)$ are uniquely determined up to multiplicative constant with
their zeros by Hadamard's factorization theorem [46]

$\Delta\left(  \lambda\right)  =C%
{\displaystyle\prod\limits_{n=-\infty}^{\infty}}
\left(  1-\dfrac{\lambda}{\lambda_{n}}\right)  ,$

$\Delta_{1}(\tau)=C_{1}%
{\displaystyle\prod\limits_{n=-\infty}^{\infty}}
\left(  1-\dfrac{\tau}{\tau_{n}}\right)  ,$\newline where $C$ and $C_{1}$ are
constants depend on $\left\{  \lambda_{n}\right\}  $ and $\left\{  \tau
_{n}\right\}  $, respectively. When$\smallskip$\newline$\lambda_{n}%
=\tilde{\lambda}_{n}$ and $\tau_{n}=\tilde{\tau}_{n}$ for all $n$,
$\Delta(\lambda)\equiv\tilde{\Delta}(\lambda)$ and $\Delta_{1}(\tau
)\equiv\tilde{\Delta}_{1}(\tau)$. Hence,$\smallskip$ $\psi_{n1}(a,\tau
)=\tilde{\psi}_{n1}(a,\tau)$. As a result, we get $M(\lambda)=\tilde
{M}(\lambda)\smallskip$. So, the proof is completed by Theorem 3.\bigskip

\begin{center}
\textbf{REFERENCES}\bigskip
\end{center}

[1] Levitan B. M. and Sargsyan I. S., Sturm-Liouville and Dirac Operators [in
Russian], Nauka (Moscow,1988).

[2] Berezanskii Yu. M., "Uniqueness theorem in the inverse spectral problem
for the Schr\"{o}dinger equation", Tr. Mosk. Mat. Obshch., 7, 3-51 (1958).

[3] Gasymov M. G. and Dzhabiev T. T., Determination of a system of Dirac
differential equations using two spectra, in Proceedings of School-Seminar on
the Spectral Theory of Operators and Representations of Group Theory [in
Russian] (Elm, Baku, 1975), pp. 46-71.

[4] Marchenko V. A., Sturm-Liouville Operators and Their Applications [in
Russian], Naukova Dumka, Kiev (1977).

[5] Nizhnik L. P., Inverse Scattering Problems for Hyperbolic Equations [in
Russian], Naukova Dumka, Kiev (1977).

[6] Gasymov M. G., Inverse problem of the scattering theory for Dirac system
of order 2n, Tr. Mosk. Mat. Obshch., 19 (1968), 41-112; Birkhauser (Basel, 1997).

[7] Guseinov I. M., On the representation of Jost solutions of a system of
Dirac differential equations with discontinuous coefficients, Izv. Akad. Nauk
Azerb. SSR, 5 (1999), 41-45.

[8] Fulton C. T., Two-point boundary value problems with eigenvalue parameter
contained in the boundary conditions, Proc. Roy. Soc. Edin. 77A, pp. 293-308, 1977.

[9] Shkalikov A. A., Boundary Value Problems For Ordinary Differential
Equations with a Parameter in Boundary Conditions Trudy Sem. Imeny I. G.
Petrovskogo, 9, pp. 190-229, 1983.

[10] Yakubov S, Completeness of Root Functions of Regular Differential
Operators, Logman, Scientific and Technical, New York, 1994.

[11] Kerimov NB, Memedov KhK, On a boundary value problem with a spectral
parameter in the boundary conditions, Sibir. Matem. Zhurnal., 40 (2), pp.
325-335, 1999 English translation: Siberian Mathematical Journal , 40 (2), pp.
281-290, 1999.

[12] Binding P. A., Browne P. J. and Watson B. A., Sturm-Liouville problems
with boundary conditions rationally dependent on the eigenparameter II.
Journal of Computational and Applied Mathematics, 148, pp. 147-169, 2002.

[13] Mukhtarov OSh, Kadakal M, Muhtarov FS. On discontinuous Sturm-Liouville
problem with transmission conditions, Journal of Mathematics of Kyoto
University, 444, pp. 779-798, 2004.

[14] Tun\c{c} E, Muhtarov OSh. Fundamental solution and eigenvalues of one
boundary value problem with transmission conditions, Applied Mathematics and
Computation , 157, pp. 347-355, 2004.

[15] Akdo\u{g}an Z, Demirci M, Mukhtarov OSh, Sturm-Liouville problems with
eigendependent boundary and transmissions conditions, Acta Mathematica
Scientia, 25B (4) pp. 731-740, 2005.

[16] Akdo\u{g}an Z, Demirci M, Mukhtarov OSh, Discontinuous Sturm-Liouville
problem with eigenparameter-dependent boundary and transmission conditions,
Acta Applicandae Mathematicae, 86 pp. 329-344, 2005.

[17] Fulton C. T.,Singular eigenvalue problems with eigenvalue parameter
contained in the boundary conditions, Proc. R. Soc. Edinb. A. 87, pp. 1-34, 1980.

[18] Amirov R. Kh., Ozkan A. S. and Keskin B., Inverse problems for impulsive
Sturm-Liouville operator with spectral parameter linearly contained in
boundary conditions, Integral Transforms and Special Functions, 20 (2009), 607-618.

[19] Guliyev N. J., Inverse eigenvalue problems for Sturm-Liouville equations
with spectral parameter linearly contained in one of the boundary conditions,
Inverse Problems, 21 (2005), 1315-1330.

[20] Mukhtarov O. Sh., Discontinuous boundary value problem with spectral
parameter in boundary conditions, Turkish J. Math., 18 (1994), 183-192.

[21] Russakovskii E. M., Operator treatment of boundary problems with spectral
parameters entering via polynomials in the boundary conditions, Funct. Analy.
Appl. 9, pp. 358-359, 1975.

[22] Binding P. A., Browne P. J. and Seddighi K., Sturm-Liouville problems
with eigenparameter dependent boundary conditions, Proc. Edinburgh Math. Soc.,
(2), 37 (1993), 57-72.

[23] Russakovskii E. M., Matrix boundary value problems with eigenvalue
dependent boundary conditions, Oper. Theory Adv. Appl. (1995).

[24] Mennicken R., Schmid H. and Shkalikov A. A., On the eigenvalue
accumulation of Sturm-Liouville problems depending nonlinearly on the spectral
parameter, Math. Nachr., 189 (1998), 157-170.

[25] Binding P. A., Browne P. J. and Watson B. A., Inverse spectral problems
for Sturm-Liouville equations with eigenparameter dependent boundary
conditions, J. London Math. Soc., 62 (2000), 161-182.

[26] Schmid H. and Tretter C., Singular Dirac systems and Sturm-Liouville
problems nonlinear in the spectral parameter, J. Differen. Equations 181 (2)
pp. 511-542, 2002.

[27] Binding P. A., Browne P. J. and Watson B. A., Equivalence of inverse
Sturm-Liouville problems with boundary conditions rationally dependent on the
eigenparameter, J. Math. Anal. Appl., 29 (2004), 246-261.

[28] Hald O. H., Discontinuous inverse eigenvalue problems, Comm. Pure Appl.
Math., 37, 539-577 (1984).

[29] Kobayashi, M., A uniqueness proof for discontinuous inverse
Sturm-Liouville problems with symmetric potentials, Inverse Problem, 5, No. 5,
767-781 (1989).

[30] Shepelsky D., The inverse problem of reconstruction of the medium's
conductivity in a class of discontinuous and increasing functions, Spectral
Oper. Theory Rel. Topics: Adv. Sov. Math., 19, 209-232 (1994).

[31] Amirov R. Kh., G\"{u}ld\"{u} Y., \.{I}nverse Problems For Dirac Operator
With Discontinuity Conditions Inside An Interval, International Journal of
Pure and Applied Mathematics, 37, No. 2, 2007, 215-226.

[32] Likov A. V. and Mikhailov Yu. A., The Theory of Heat and Mass Transfer
Qosenergaizdat, 1963 (Russian).

[33] Meschanov V.P. and Feldstein A. L., Automatic Design of Directional
Couplers, Sviaz, Moscow, 1980.

[34] Tikhonov A. N. and Samarskii A. A., Equations of Mathematical Physics,
Pergamon, Oxford, 1990.

[35] McLaughlin J. and Polyakov P., On the uniqueness of a spherical symmetric
speed of sound from transmission eigenvalues, J. Differ. Eqns 107 (1994), pp. 351-382.

[36] Voitovich N. N., Katsenelbaum B. Z. and Sivov A. N., Generalized Method
of Eigen-vibration in the Theory of Diffraction, Nauka, Moskov, 1997 (in Russian).

[37] Titeux I, Yakubov Ya., Completeness of root functions for thermal
conduction in a strip with peicewise continuous coefficients, Mathematical
Models and Methods in Applied Sciences 1997; 7 1035-1050.

[38] Yurko V. A., Integral transforms connected with discontinuous boundary
value problems, Integral Transforms Spec. Funct. 10 (2000), pp. 141-164.

[39] Freiling G. and Yurko V. A., Inverse Sturm-Liouville Problems and Their
Applications, Nova Science, New York, 2001.

[40] Kadakal M., Mukhtarov O. Sh., Sturm-Liouville problems with
discontinuities at two points, Comput. Math. Appl. 54 (2007) 1367-1379.

[41] Yang Qiuxia and Wang Wanyi , Asymptotic behavior of a differential
operator with discontinuities at two points, Mathematical Methods in the
Applied Sciences, 34 pp. 373-383, 2011.

[42] Shahriari Mohammad , Akbarfam Aliasghar Jodayree , Teschl Gerald ,
Uniqueness for inverse Sturm-Liouville problems with a finite number of
transmission conditions, J. Math. Anal. Appl. 395 (2012) 19-29.

[43] G\"{u}ld\"{u} Y., Inverse eigenvalue problems for a discontinuous
Sturm-Liouville operator with two discontinuities, Boundary Value Problems
2013, 2013:209.

[44] Yang Chuan-Fu , Uniqueness theorems for differential pencils with
eigenparameter boundary conditions and transmission conditions, J.
Differential Equations 255 (2013) 2615--2635.

[45] Zhdanovich V. F., Formulae for the zeros of Dirichlet polynomials and
quasi-polynomials, Dokl. Acad. Nauk SSSR 135 (8) (1960), pp. 1046-1049.

[46] Titchmarsh E. C., The Theory of Functions, Oxford University Press,
London (1939).

\end{document}